\documentclass{article}
\usepackage[final]{neurips_2022}%

\usepackage[utf8]{inputenc} %
\usepackage[T1]{fontenc}    %
\usepackage{hyperref}       %
\usepackage{url}            %
\usepackage{booktabs}       %
\usepackage{nicefrac}       %
\usepackage{microtype}      %
\usepackage{graphicx}
\usepackage{subfigure}
\usepackage{hyperref}
\usepackage{times}		 
\DeclareUnicodeCharacter{00A0}{~}
\usepackage{amssymb,amsmath,amscd,amsfonts,amsthm,bbm,mathrsfs,yhmath}
\usepackage[shortlabels]{enumitem}
\usepackage{hyperref}
\usepackage{caption}
\usepackage{graphics}
\usepackage{mathtools}
\usepackage{cleveref}
\usepackage{macros}

\usepackage{todonotes}

\theoremstyle{definition}
\newtheorem{theorem}{Theorem}
\newtheorem{lemma}[theorem]{Lemma}

\newtheorem{proposition}[theorem]{Proposition}
\newtheorem{definition}{Definition}
\newtheorem{remark}{Remark}

\newtheorem{example}{Example}
\newenvironment{sproof}{%
  \proof}{\endproof}

\newcommand{\subscript}[2]{$#1 _ #2$}
\newlist{assumplist}{enumerate}{1}
\setlist[assumplist]{label=(\subscript{\textbf{A}}{{\arabic*}})}
\Crefname{assumplisti}{Assumption}{Assumptions}

\newlist{assumplist2}{enumerate}{1}
\setlist[assumplist2]{label=(\subscript{\textbf{A'}}{{\arabic*}})}
\Crefname{assumplist2i}{Assumption}{Assumptions}

\newlist{assumplist3}{enumerate}{1}
\setlist[assumplist3]{label=(\subscript{\textbf{A''}}{{\arabic*}})}
\Crefname{assumplist3i}{Assumption}{Assumptions}

\title{Mirror Descent with Relative Smoothness in Measure Spaces, with application to Sinkhorn and EM}
\author{
 Pierre-Cyril Aubin-Frankowski\\ DI ENS, Ecole normale supérieure,\\ Université PSL, CNRS, INRIA Paris\\ \texttt{pierre-cyril.aubin@inria.fr} \And Anna Korba\\ CREST, ENSAE\\ IP Paris\\ \texttt{anna.korba@ensae.fr} \And Flavien Léger \\ INRIA Paris \\ \texttt{flavien.leger@inria.fr } 
}
\begin{document}

	\maketitle
	
	\begin{abstract}
	    Many problems in machine learning can be formulated as optimizing a convex functional over a vector space of measures. This paper studies the convergence of the mirror descent algorithm in this infinite-dimensional setting. Defining Bregman divergences through directional derivatives, we derive the convergence of the scheme for relatively smooth and convex pairs of functionals. Such assumptions allow to handle non-smooth functionals such as the Kullback--Leibler (KL) divergence. Applying our result to joint distributions and KL, we show that Sinkhorn's primal iterations for entropic optimal transport in the continuous setting correspond to a mirror descent, and we obtain a new proof of its (sub)linear convergence. We also show that Expectation Maximization (EM) can always formally be written as a mirror descent. When optimizing only on the latent distribution while fixing the mixtures parameters -- which corresponds to the Richardson--Lucy deconvolution scheme in signal processing -- we derive sublinear rates of convergence.%

	\end{abstract}

\section{Introduction}
Many important problems in machine learning and computational statistics can be cast as an optimization problem over the space of probability distributions, where the
objective functional assesses the dissimilarity to
a target distribution $\tg$ on $\R^d$. Classical dissimilarities include $f$-divergences, Integral Probability Metrics (IPMs), or optimal transport distances among others.  In Bayesian inference for instance, it is common to optimize the Kullback--Leibler (KL) divergence to the target, which corresponds to the posterior distribution of the parameters of interest. 
In  generative modelling, the goal is to generate data whose distribution is similar to the training set distribution defined by samples of the target, where the similarity is often measured by an integral probability metric or an optimal transport distance \citep{arjovsky2017wasserstein,dziugaite2015training}.  
In supervised learning, optimizing an infinite-width one hidden layer neural network, e.g. through the mean squared error, corresponds to minimizing a functional on the space of probability distributions over the parameters of the network \citep{chizat2018global,mei2018mean,rotskoff2018neural}. In particular, the objective functional can be identified to a Maximum Mean Discrepancy (MMD) in the well-specified setting \citep{arbel2019maximum}. Many other problems in machine learning can be formalized in this framework \citep{chu2019proba}.

Once the objective functional is chosen, one has to select an optimization algorithm that is well-suited to the geometry of the problem. In this article, we consider the widely used mirror descent scheme, a first-order optimization method based on Bregman divergences. While the traditional smoothness and strong convexity assumptions, required by standard convergence analysis, do not always hold over measure spaces, their ``relative'' versions have received increased interest. The relative smoothness assumption was first suggested by \citet{Birnbaum2011} in the context of algorithmic game theory, but remained unnoticed by the optimization community, until \citet{bauschke2017descent} discovered the same concept independently, while \citet{lu2018relatively} coupled it with relative strong convexity. We extend their work to the infinite dimensional setting and target specifically the KL divergence. Indeed, when using the entropy as Bregman divergence, mirror descent is known to yield multiplicative updates for the measures. In this paper we study two such schemes in machine learning, i.e. Sinkhorn's algorithm, widely used to solve entropic optimal transport \citep{peyre2019computational}, and the EM algorithm, a very common approach to fitting probabilistic models.

\textbf{Related work.} \cite{chizat2021convergence} gave convergence rates of mirror descent on measure spaces for integral functionals with Lipschitz gradients, working mostly in the $L^1$ space and without leveraging relative smoothness. However his assumptions do not cover $f$-divergences, such as the ubiquitous $\KL$, nor entropic regularized transport as in \cite{Leger2020}. His setting also implies Gâteaux differentiability, which is classical for mirror descent in Banach spaces \citep{Bauschke2003}. As the entropy is not differentiable in infinite dimensions, directional derivatives were instead used by \citet{resmerita2005regularization} for mirror descent, and applied by \cite{chu2019proba} to gradient descent.  Connections between mirror descent and Sinkhorn iterations for the entropic optimal transport problem were first investigated in \cite{mishchenko2019sinkhorn,menschpeyre, Leger2020}. Our framework is closest to \citet[]{Leger2020}, which we simplify by considering primal rather than dual iterations, and extend by also deriving a linear convergence rate. \citet[]{Kunstner2021Homeomorphic} recently showed that EM over parametric exponential mixtures corresponds to a mirror descent scheme with relative smoothness properties. Their setting is complementary to ours since we consider instead fixed mixtures and a nonparametric latent distribution. A prominent alternative setting to ours when optimizing over measures is based on (grid-free) Wasserstein gradient flows, which we do not cover here (see \Cref{sec:related_work} for a discussion on the difference of geometries). A remarkable discussion on the optimization over measures with different geometries can be found in \citet{trillos2020bayesian}.

\textbf{Contributions.} We propose a rigourous framework for the analysis of the infinite-dimensional version of mirror descent over measure spaces. In this setting, we recover the rate of convergence of mirror descent under relative smoothness and convexity, previously shown in finite dimensions. Defining Bregman divergences through directional derivatives, we are able to consider objective functionals over measures that are not smooth in the "standard" sense, but satisfy relative smoothness and/or convexity with respect to a Bregman divergence, e.g.\  KL. Focusing on optimization over joint distributions, we show that both Sinkhorn's primal iterations for entropic optimal transport in the continuous setting and EM can be written as a mirror descent. We then obtain a new proof of Sinkhorn's (sub)linear convergence. For EM, when optimizing on the latent distribution while fixing the mixtures, a choice which coincides with Richardson--Lucy deconvolution, we derive new sublinear rates of convergence.

This paper is organized as follows. \Cref{sec:background} introduces the necessary background on derivatives in measure spaces and relative smoothness and convexity. \Cref{sec:MD_convergence} discusses the well-posedness of the mirror descent scheme and provides our proof of convergence adapted from \cite{lu2018relatively}. \Cref{sec:MD_algorithms} recovers the convergence of algorithms such as Sinkhorn's iterations and latent EM as special cases of mirror descent with relative smoothness and convexity. %

	\section{Background and definitions}\label{sec:background}
	
	In this section, we set the mathematical framework in which we will rigourously reformulate relative smoothness and convexity on a space of measures.%
	
	\textbf{Notation.} Given a topological vector space $\Y$ with topology $\tau$, the domain $\dom(f)$ of an extended-valued function $f:\Y\rightarrow\R\cup\{\pm\infty\}$ is the set of points of $\Y$ where $f$ takes finite values. The function $f$ is said to be proper if $\dom(f)$ is non-empty and if $f$ never takes the value $-\infty$. It is $\tau$-lower semicontinuous (l.s.c.) if its sublevel sets are $\tau$-closed. We consider a dual pair $(\Y,\Y^*)$ with duality product $\scalidx{\cdot}{\cdot}{\Y^*\times \Y}$ which induces a $\Y^*$-weak topology on $\Y$ \citep[see][]{aliprantis2006infinite}. We write $\Int C$ for the interior of a set $C\subset \Y$. Let $\X \subset \R^d$, and fix a locally convex topological vector space of measures $(\cM(\X),\tau)$, which could be for instance $L^1(\d\rho)$, $L^2(\d\rho)$ where $\rho$ is a reference measure, or the space of Radon measures $\cM_r(\X)$ with the total variation (TV) norm. Fix $\cM^*(\X)$ a topological dual of $\cM(\X)$. For $\mu\in \cM(\X)$ and $f\in \cM^*(\X)$, we use the shorthand $\ps{f, \mu} =\ps{f, \mu}_{\cM^*(\X)\times \cM(\X)}$, formally equal to $\int_{\X}f(x)\mu(dx)$. We denote by $\cP(\X)$  the subset of measures $\mu\in \cM(X)$ with mass 1, and, for any $\mu,\nu\in \cM(\X)$, we write $\mu \ll \nu$ when $\mu$ is absolutely continuous w.r.t\ $\nu$, i.e.\ when it has a Radon--Nikodym derivative $d\mu/d\nu$.

	Consider a convex functional $\cF:\cM(\X)\rightarrow \R\cup \{+\infty\}$ and the following minimization problem
	\begin{equation}\label{eq:opt-cons}
		\begin{aligned}
			\min_{\nu \in C} \cF(\nu)
		\end{aligned}
	\end{equation}
	where $C\subset \cM(\X)$ is a convex set. To solve this convex optimization problem, a classical choice is to resort to a mirror descent scheme \citep[see e.g.][]{beck2003mirror}. The latter is a first-order optimization scheme based on the knowledge of the ``derivative'' of the objective functional $\cF$ at each iteration. The difficulty is to choose the appropriate notion of derivative. Indeed, Gâteaux and Fréchet derivatives have to be defined in every direction (see~\Cref{sec:add_defs}), thus requiring that the points of differentiability belong to the interior of the domain $\Int(\dom(\cF))$ of the functional $\cF$ considered. In infinite dimensions, for functionals defined on positive measures, such as the negative entropy, $\Int(\dom(\cF))$ is however empty\footnote{ Intuitively if $\X$ contains an open set, then any positive measure can be infinitesimally perturbed to have negative values. For finite sets $\X$, this phenomenon does not occur, see also \Cref{rmk:KL_properties}.}. Consequently, following \cite{resmerita2005regularization}, we favor a weaker notion, that of directional derivatives.\footnote{In finite dimensions, \citet{Maddison2021} also defined Bregman divergences through directional derivatives, but under the stringent assumption of essentially smooth convex functions, an assumption which does not extend well to infinite dimensions.} This comes at the price of manipulating $\pm\infty$ values but ensures that the considered derivatives are always well-defined for convex functionals. Besides, whenever the directional derivative is a linear function in a restricted set of directions, the notion of first variation that comes next will enable us to perform all the computations we need, as if the function was Gâteaux differentiable.%

	\begin{definition}[Directional derivative]\label{def:directional_derivative}
		If it exists, the \emph{directional derivative} of  $\cF:\cM(\X) \rightarrow\R\cup\{\pm\infty\}$ at a point $\nu\in \dom(\cF)$ in the direction $\mu\in\cM(\X)$ is defined as
		\begin{equation}
			d^+\!\cF(\nu)(\mu)=\lim_{h \rightarrow 0^+}\frac{\cF(\nu+h \mu)- \cF(\nu)}{h}.
		\end{equation}
	\end{definition}
	
	\begin{remark}%
	\label{rmk:directional_derivative}
		In particular, for convex and proper functions, $d^+\!\cF(\nu)(\mu)$ exists and belongs to $\R\cup\{\pm\infty\}$ \cite[Lemma 7.14]{aliprantis2006infinite}. This is a direct consequence of the nondecreasingness of $\R_+^* \ni h \mapsto\frac{\cF(\nu+h \mu)- \cF(\nu)}{h}$ for any $\nu,\mu$ and convex $\cF$. The monotonicity also entails that $d^+\!\cF(\nu)(\mu)\le \cF(\nu+\mu)- \cF(\nu)<\infty$ whenever $\nu$ and $\nu+\mu$ belong to $\dom(\cF)$ and that $d^+$ is a linear operation over the cone of convex functions. Note that $\mu\mapsto d^+\!\cF(\nu)(\mu)$ is not always $\tau$-l.s.c.\ although it is always positively homogeneous, and, whenever $\cF$ is convex, it is convex. %

	\end{remark}	
    
	\begin{definition}[First variation]\label{def:first_var_dual}
		Let $\cF:\cM(\X)\rightarrow \R\cup \{+\infty\}$ be a functional and $C$ be a subset of $\cM(\X)$. If it exists, the \emph{first variation} of $\cF$ over $C$ evaluated at $\mu\in \dom(\cF)\cap C$ is the element $\nabla_{\!C} \cF(\mu)\in \cM^*(\X)$, unique up to orthogonal components to $\Sp(\dom(\cF)\cap C-\mu)$, such that:
		\begin{equation}\label{eq:first_var_dual}
		\ps{\nabla_{\!C} \cF(\mu),\xi}
			=d^+\!\cF(\mu)(\xi)
		\end{equation} 
		for all $\xi= \nu-\mu \in \cM(\X)$, where $\nu \in \dom(\cF)\cap C$.
	\end{definition}

	By \Cref{rmk:directional_derivative}, we have that $d^+\!\cF(\mu)(\xi)\in[-\infty,\infty)$%
	, since $\nu,\mu\in \dom(\cF)\cap C$. Naturally, if $\cF$ has a Fréchet or Gâteaux derivative at $\mu$, which implies that $\mu\in \Int(\dom(\cF))$, then it coincides with the first variation at $\mu$. Calling first variations the derivatives of functionals stems from the field of calculus of variations; a pragmatic approach when dealing with probability measures can be found in \citet[Definition 7.12]{Santambrogio2015}, where $\nabla_{\cP(\X)} \cF(\mu)$ is defined as a measurable function, we instead take it most in often in $L^\infty(\X)$. We now introduce Bregman divergences over measures through directional derivatives.
	\begin{definition}\label{def:bregman_div}(Bregman divergence)
		Let $\phi:\cM(\X)\rightarrow \R\cup\{+\infty\}$ be a convex functional. For $\mu \in \dom(\phi)$, the $\phi$-\emph{Bregman divergence} is defined for all $\nu\in\dom(\phi)$ by
		\begin{equation}\label{eq:div_Bregman}
			D_{\phi}(\nu|\mu)=\phi(\nu)-\phi(\mu)-d^+\phi(\mu)(\nu- \mu)\in [0,+\infty],
		\end{equation}
		and $+\infty$ elsewhere. The function $\phi$ is referred to as \textit{the Bregman potential}.
	\end{definition}%
	
	In the previous definition, the restriction to $\dom(\phi)$ is necessary to avoid the substraction of infinite values. As a direct consequence of \Cref{lem:conv_directionnal_deriv} in Appendix, a stricly convex $\phi$ entails that the Bregman divergence $D_{\phi}$ separates measures, i.e.\ $D_{\phi}(\nu|\mu)=0$ if and only if $\nu=\mu$. Note that $D_{\phi}(\cdot |\mu)$ is a difference of convex functions, so it is not convex in general. Nevertheless the existence of a first variation \eqref{eq:first_var_dual} of $\phi$ over $C$, resulting in the last term in \eqref{eq:div_Bregman} being linear, is sufficient to ensure the convexity of the restriction of $D_{\phi}(\cdot |\mu)$ to $C$. Bregman divergences have useful immediate properties: since $d^+$ is a linear operation over convex functions, so is the Bregman divergence, i.e.\ for two convex $\phi,\psi$,  $D_{\phi+\psi}=D_{\phi}+D_{\psi}$. Moreover, it is idempotent, as shown in the following lemma.
		\begin{lemma}[Idempotence of Bregman divergence]\label{lem:Bregman_translated} %
	Let $\phi:\cM(\X)\rightarrow \R\cup\{+\infty\}$ be a convex functional. Assume that given $\xi\in\dom(\phi)$, the first variation $ \nabla_{\!C} \phi(\xi)$ exists, then, for all $\mu,\nu\in C\cap \dom(\phi)$, %
	    $D_{D_{\phi}(\cdot|\xi)}(\nu|\mu)=D_{\phi}(\nu|\mu)$.
    \end{lemma}
    \begin{proof}  %
    Since $\psi:\tilde{\mu}\mapsto -d^+\phi(\xi)(\tilde{\mu}-\xi)=-
    \ps{\nabla_{\!C} \phi(\xi),\tilde{\mu}-\xi}
    $ is convex over $C\cap\dom(\phi)$, we can apply the linearity of the Bregman divergence: 
        \begin{equation*}
	    D_{D_{\phi}(\cdot|\xi)}(\nu|\mu)=D_{\phi}(\nu|\mu)+D_{-\phi(\xi)}(\nu|\mu)+D_{\psi}(\nu|\mu)=D_{\phi}(\nu|\mu),%
	    \end{equation*}
	   since the Bregman divergence of a constant or of a linear form is null.%
    \end{proof}

We are now ready to introduce the notions of relative smoothness \citep{bauschke2017descent} and convexity \citep{lu2018relatively} of a functional w.r.t\ a Bregman potential.
	\begin{definition} \label{def:relative_smoothness}(Relative smoothness and convexity)
	Let $\cF:\cM(\X)\to \R\cup\{+\infty\}$ be a convex proper functional.  Given a scalar $L\ge 0$, we say that $\cF$ is $L$-smooth relative to $\phi$ over $C$ if, for any $\mu, \nu \in \dom(\cF)\cap\dom(\phi)\cap C$, we have
		\begin{equation}\label{eq:relative_smoothness}%
			D_{\cF}(\nu| \mu)=\cF(\nu) - \cF(\mu) - d^+\!\cF(\mu)(\nu- \mu) \le L D_{\phi}(\nu| \mu).
		\end{equation}
	Conversely, we say that $\cF$ is $l$-strongly convex relative to $\phi$ over $C$, for some scalar $l\ge 0$, if, for any $\mu, \nu \in \dom(\cF)\cap\dom(\phi)\cap C$, we have
		\begin{equation}\label{eq:relative_convexity}
			D_{\cF}(\nu| \mu) \ge l D_{\phi}(\nu| \mu).
		\end{equation}
\end{definition}
\begin{example}[L-smoothness]\label{ex:l_smoothness}
Choosing $\phi(\mu)=\|\mu\|_{\cM(\X)}^2$ the square norm on $\cM(\X)$ shows that relative smoothness %
extends the notion of $L$-smooth functionals (see for instance \citet[Lemma 3.1]{chizat2021convergence}), i.e.\ functionals with $L$-Lipschitz Gâteaux derivative, thus satisfying:
\begin{equation}\label{eq:smoothness}
    \cF(\nu) -  \cF(\mu) -d^+\cF(\mu)(\nu- \mu) \le L \|\nu - \mu\|^2.
\end{equation}
\end{example}

Notice that by \Cref{lem:Bregman_translated}, provided $ \nabla_{\!C} \phi(\xi)$ is well-defined, a Bregman divergence objective $D_{\phi}(\cdot |\xi)$ is always 1-relatively smooth and strongly convex w.r.t.\ $\phi$. This we will heavily exploit for mirror descent schemes that involve the KL divergence both as an objective and Bregman divergence in \Cref{sec:MD_algorithms}. Interestingly, relative smoothness and convexity can be characterized in different, equivalent ways, see \Cref{lem:equivalences} in the Appendix. We now turn to the analysis of the mirror descent scheme using the above framework. %

\section{Mirror descent over measures and convergence}\label{sec:MD_convergence}

In the following, $\phi$ is assumed to be strictly convex.
The relative smoothness assumption \eqref{eq:relative_smoothness} of the convex objective functional $\cF$ w.r.t.\ a Bregman potential $\phi$ suggests to minimize iteratively over $\nu\in C$ the function $\nu \mapsto \cF(\mu) + d^+\!\cF(\mu)(\nu- \mu) + L D_{\phi}(\nu| \mu)$, acting as an upper approximation of $\cF(\nu)$. Starting from a given $\mu_0 \in \cM(\X)$, the mirror descent iterates are thus defined at each time $n\ge 0$ as
\begin{equation}\label{eq:Bregman_prox_algo}
	\hspace*{-0.3cm} \mu_{n+1}=\argmin_{\nu \in C}\{ 
	d^+\!\cF(\mu_n)(\nu- \mu_n)
	+ L D_{\phi}(\nu|\mu_n) \}.
\end{equation}
Let $\cR\subset C$ be a given subset. %
As proven later in this section, sufficient conditions for the convergence of the scheme \eqref{eq:Bregman_prox_algo} are:
\begin{assumplist}[leftmargin=7mm]
	\setlength\itemsep{0.2em}
	\item \label{ass:diff_F_phi}(Existence) The sequence of iterates $(\mu_{n})_{n\in\N}$ defined by \eqref{eq:Bregman_prox_algo} exist, belong to $\cR$, and are unique.
	\item \label{ass:relative_smoothness}(Relative smoothness/convexity) For some $l,L\ge 0$, the functional $\cF$ is $L$-smooth and $l$-strongly convex relative to $\phi$ as in \Cref{def:relative_smoothness} for elements of $\cR$.
	\item \label{ass:first-variation}(Existence of first variation of $\phi$) For each $n\ge0$, the first variation $\nabla_{\!C}\phi(\mu_n)$ exists.
\end{assumplist}
These assumptions have to be verified on a case-by-case basis. In the simplest case, one can take $\cR=C$. However the set $\cR$ does not have to be convex (see \Cref{sec:MD_algorithms}). On the other hand,  \Cref{ass:diff+inf-compact} below, ensures that the iterates in \eqref{eq:Bregman_prox_algo} are well-defined, uniqueness resulting from the strict convexity of $\phi$. 
\begin{assumplist2}[leftmargin=8mm]
	\setlength\itemsep{0.2em}
	\item \label{ass:diff+inf-compact}(Lower semicontinuity and coercivity)  (i) the set $C$ is $\tau$-closed in $\cM(\X)$, the functionals $\cG_n(\cdot):= d^+\!\cF(\mu_n)(\cdot- \mu_n)$ and $D_{\phi}(\cdot | \mu_n)$ are proper and $\tau$-l.s.c.\ when restricted to $C$, and the functional $\cG_n+D_{\phi}(\cdot |\mu_n)+i_{C}$\footnote{$i_{C}$ denotes the indicator function of the set $C$, defined by $i_{C}(\mu)=0$ if $\mu\in C$, $+\infty$ otherwise for any $\mu\in \cM(\X)$. Notice that $i_{C}$ being $\tau$-l.s.c.\ is equivalent to $C$ being $\tau$-closed in $\cM(\X).$} has at least one $\tau$-compact sublevel set. (ii) For each $n\ge0$, the first variations $\nabla_{\!C}\phi(\mu_n)$ exist. (iii) The iterates belong to $\cR$.
\end{assumplist2}

By \citet[Theorem 3.2.2]{attouch2014variational}, \Cref{ass:diff+inf-compact} implies  \ref{ass:diff_F_phi}. Indeed, (i) $\tau$-lower semicontinuity and $\tau$-compactness guarantee the existence of minimizers of the objective \eqref{eq:Bregman_prox_algo}, while (ii) the existence of $\nabla_{\!C}\phi(\mu_n)$ guarantees that $D_{\phi}(\cdot |\mu_n)$ is strictly convex (since $\phi$ is assumed strictly convex) hence unicity of the minimizer. Regarding \Cref{ass:diff+inf-compact}(i), notice that $\cG_n$ is proper and $\cM^*(\X)$-weak-l.s.c.\ as soon as $\cF$ has a first variation at $\mu_n$, since in this case $\cG_n$ is linear on $C$. We refer to \Cref{sec:assumptions_details} for more details on proving that \Cref{ass:diff+inf-compact} holds for some $\cF$, $\phi$ and $C$. %

Mirror descent can also be defined through a subdifferential constraint if $\cF$ also has first variations.
\begin{lemma}[Mirror descent dual iteration]\label{rmk:dual_iteration} If $\nabla_{\!C}\cF(\mu_n)$ and $\nabla_{\!C}\phi(\mu_n)$ exist for all $n\ge 0$, then \eqref{eq:Bregman_prox_algo} is equivalent to
\begin{equation}\label{eq:dual_iteration}
\nabla_{\!C} \phi(\mu_{n})-\frac{1}{L}\nabla_{\!C}\cF(\mu_n) \in \partial_{\!C} \phi(\mu_{n+1}):=\{p \, | \,\forall \nu\in C,\, \ps{p,\nu-\mu_{n+1}}\le d^+\!\phi(\mu_{n+1})(\nu- \mu_{n+1}) \}
\end{equation}
Thus, if $\partial_{\!C} \phi(\mu_{n+1})=\{\nabla_{\!C} \phi(\mu_{n+1})\}$,  \eqref{eq:Bregman_prox_algo} corresponds to $\nabla_{\!C} \phi(\mu_{n+1}) - \nabla_{\!C} \phi(\mu_n) =-\frac{1}{L}\nabla\cF(\mu_n)$.
\end{lemma}
\begin{proof}
    The minimization \eqref{eq:Bregman_prox_algo} is equivalent to having, for all $\nu\in C$,
    \begin{gather*}
        d^+\!\cF(\mu_n)(\nu- \mu_n)
	+ L D_{\phi}(\nu|\mu_n) \ge d^+\!\cF(\mu_n)(\mu_{n+1}- \mu_n)
	+ L D_{\phi}(\mu_{n+1}|\mu_n)\\
	\ps{\nabla_{\!C}\cF(\mu_n) -L\nabla_{\!C} \phi(\mu_{n}),\nu-\mu_{n+1}}+L(\phi(\nu)-\phi(\mu_{n+1}))\ge 0.
    \end{gather*}
    Take $\tilde \nu\in C$, set $\nu=\mu_{n+1}+t(\tilde \nu -\mu_{n+1})$ for $t\in[0,1]$. Taking the limit $t\rightarrow 0^+$ yields the result.
\end{proof}
    In the general case, as discussed in finite dimensions in \citet[Remark 3]{bauschke2017descent}, one needs extra assumptions to justify that $\mu_{n+1}$ exists in \eqref{eq:dual_iteration}, akin to the invertibility of $\nabla \phi$ or that $\phi$ is essentially smooth or of Legendre type. To avoid any restrictive assumption required to use \eqref{eq:dual_iteration}, we stick with the minimal formulation \eqref{eq:Bregman_prox_algo} as was also done by \citet[]{bauschke2017descent,lu2018relatively}.

We now state a preliminary result, %
known as the "three-point inequality" or "Bregman proximal inequality" in the optimization literature \cite[Lemma 3.2]{chen1993convergence}, \cite[Lemma 1]{lan2011PrimaldualFM},  %
useful to 
prove the convergence of the mirror descent scheme, similarly to \citet{lu2018relatively}.

\begin{lemma}[Three-point inequality]\label{lem:three-point}
	Given $\mu\in\cM(\X)$ and some proper convex functional $\cG:\cM(\X)\rightarrow \R\cup\{+\infty\}$, if $\nabla_{\!C}\phi(\mu)$ exists, as well as $\bar{\nu}=\argmin_{\nu \in C}\{\cG(\nu)+D_{\phi}(\nu | \mu)\}$, then for all $\nu\in C\cap \dom(\phi) \cap \dom(\cG)$:%
	\begin{equation}\label{eq:three_point_ineq}
		\hspace*{-0.2cm} \cG(\nu)+ D_{\phi}(\nu| \mu) \ge \cG(\bar{\nu}) + D_{\phi}(\bar{\nu}| \mu) +  D_{\phi}(\nu| \bar{\nu}).
	\end{equation}
\end{lemma}
\begin{proof}
    The existence of $\nabla_{\!C}\phi(\mu)$ entails that $C\cap\dom(D_\phi(\cdot |\mu))=C\cap\dom(\phi)$. Set $f(\cdot)=\cG(\cdot)+D_{\phi}(\cdot | \mu)$. Then, by linearity of the Bregman divergence and \Cref{lem:Bregman_translated}, we obtain that, for any $\nu\in C\cap \dom(\phi) \cap \dom(\cG)$,
		\begin{equation}\label{eq:three-point_intermediate}
	    D_{f}(\nu|\bar{\nu})=D_{\cG}(\nu|\bar{\nu})+D_{D_{\phi}(\cdot | \mu)}(\nu|\bar{\nu})=D_{\cG}(\nu|\bar{\nu})+D_{\phi}(\nu|\bar{\nu})\ge D_{\phi}(\nu|\bar{\nu}).
	    \end{equation}
     By optimality of  $\bar{\nu}$, for all $\nu\in C$, $d^+(f)(\bar{\nu})(\nu-\bar{\nu})=\lim_{h \rightarrow 0^+}(f((1-h)\bar{\nu}+h \nu)- f(\bar{\nu}))/h\ge 0$; which is equivalent to $f(\nu)\ge f(\bar{\nu}) + D_f(\nu | \bar{\nu})$. We conclude using \eqref{eq:three-point_intermediate} and the definition of $f$.%
\end{proof}
The following theorem  gives the rate of  convergence of mirror descent for relatively smooth and convex pairs of functionals, and 
extends to infinite dimensions the convergence result of \citet[Theorem 3.1]{lu2018relatively}. Its proof can be found in \Cref{sec:proof_th_rate}.  %

\begin{theorem}[Convergence rate]\label{th:rate}
Assume that Assumptions \ref{ass:diff_F_phi}, \ref{ass:relative_smoothness} and \ref{ass:first-variation} hold. Consider the mirror descent scheme \eqref{eq:Bregman_prox_algo}, then for all $n\ge 0$ and all $\nu \in \dom(\cF)\cap\dom(\phi)\cap  \cR$, we have
	\begin{equation}
		\cF(\mu_n)-\cF(\nu)\le  \frac{l D_{\phi}(\nu|\mu_0)}{\left(1+\frac{l}{L-l}\right)^n -1} \le \frac{L}{n}D_{\phi}(\nu|\mu_0),
	\end{equation}
	where, in the case $l=0$, the middle expression is defined in the limit as $l \to 0^+$.
\end{theorem}

\begin{remark}[About the proof of convergence]\label{rmk:proof_cv_main}
Our proof resembles the one of \citet{lu2018relatively}, which also relies on a three-point inequality as stated in  \Cref{lem:three-point}. However, the proof of the latter inequality in finite dimensions relies on a sum of subdifferentials formula, which is hard to verify for general functionals and in particular for the KL divergeence defined below  %
(see also \Cref{rmk:proof_cv_app} in Appendix). On the contrary, by working with directional derivatives and first variations, we circumvent  most of the difficulties related to (sub)differentiability.
\end{remark}
 An important example is the one discussed below where $\phi$ is chosen to be the negative entropy $\phine$. %

\begin{example}[The KL divergence and negative entropy]
\label{rmk:KL_properties}
The Kullback--Leibler (KL) divergence and the negative entropy are defined for $\mu \ll \tg$ and $\mu \ll \rho$, writing $\mu(x)=\nicefrac{d\mu}{d\rho}(x)$, respectively as
\begin{equation}\label{eq:kl_def}
	\KL(\mu|\tg)=\int_\X \ln\left(\nicefrac{d\mu}{d\tg}(x)\right)d\mu(x)
	,\quad
	\phine(\mu)=\int_{\X} \ln (\mu(x))\mu(x)d\rho(x)=\KL(\mu|\rho)
	, 
\end{equation}
where $\rho$ is some reference finite measure on $\X$. It is straightforward to show that $\KL$ can be written as a Bregman divergence of $\phine$ if $\mu\ll\tg\ll \rho$, i.e. $D_{\phine}(\mu|\tg)=\KL(\mu|\tg)$; hence one can choose  $\phi=\phine$ for the mirror descent scheme \eqref{eq:Bregman_prox_algo}.  \Cref{ass:diff+inf-compact}, guaranteeing that the iterates \eqref{eq:Bregman_prox_algo} are well-posed, is satisfied for instance
when $\cM(\X)=L^1(\X)$ and if there exists $\kappa_0,\kappa_1>0$ such that $\kappa_0\le\nicefrac{d\mu_n}{d\tg}(x)\le\kappa_1$ almost everywhere over $\X$ for any $n$, which is the case for Sinkhorn's and EM iterates. This implies that the first variations of the negative entropy or KL belong to $L^\infty$ at $\mu_n$, see also \Cref{rmk:KL_properties_app} in Appendix for more details. Moreover, by \Cref{rmk:dual_iteration}, exponentiating the dual iteration \eqref{eq:dual_iteration} recovers the classical multiplicative scheme: $\mu_{n+1}=\mu_n e^{-\frac{1}{L}\nabla \cF(\mu_n)}, \; n\ge 0$. This observation generalizes to the Iterative Proportional Fitting Procedure, also known as Sinkhorn's algorithm.
\end{example}

$\KL$ is a strong Bregman divergence, in the sense that it dominates a wide range of objective functionals $\mathcal{F}$. Indeed, we already know from \Cref{sec:background} and the idempotence property that $\cF = \KL(\cdot|\tg)$ is 1-relatively smooth w.r.t.\ $\phine$, since for any $\mu,\nu \in \dom(\cF)$,   $D_{\cF}(\nu|\mu) = D_{\phine}(\nu|\mu)=\KL(\nu|\mu)$.
In \Cref{sec:MD_algorithms} we will extensively use this fact (sometimes applying $\KL$ to joint distributions rather than marginals, e.g. for Sinkhorn's algorithm). This is of crucial importance because $\cF=\KL(\cdot|\tg)$ is not a smooth objective in the "standard" sense (see %
\Cref{ex:l_smoothness}) - hence convergence proofs requiring the latter cannot apply -, but it is relatively smooth w.r.t.\ itself. %
Indeed the KL diverges for Dirac masses, so is unbounded over the bounded set $\cP(\X)$, and thus $\KL$ does not have subquadratic growth \eqref{eq:smoothness} w.r.t.\ any norm on measures. 
Other objective functionals can be dominated by KL, such as the Maximum Mean Discrepancy (MMD) for bounded kernels, see  \Cref{prop:rel_smooth_mmd} in \Cref{sec:app_mmd_kl}. %

\section{Applications to Sinkhorn and EM with convergence rates}\label{sec:MD_algorithms}

We now analyze the convergence of two algorithms, Sinkhorn and Expectation-Maximization (EM), by showing that they can be written as mirror descent schemes based on the $\KL$ divergence, in order to apply the results of \Cref{sec:MD_convergence}. Note that the relative smoothness was first introduced by \citet{Birnbaum2011} for this very purpose, to study the convergence of Proportional Response Dynamics.

In both the Sinkhorn and EM settings we will be given two probability spaces $(\X,\tgX)$ and $(\Y,\tgY)$. We recall that $\cP(\X\times \Y)$ denotes the subset of Radon measures $\cM_r(\X\times \Y)$ with mass 1. We equip $\cM_r$ with $L^\infty(\X\times \Y)$ as dual space. A joint measure $\pi\in\cP(\X\times \Y)$ is also called a \emph{coupling} between its first $p_\X\pi$ and second $p_\Y\pi$ marginals. We denote by $\Pi(\tgX,*)$ the set of couplings having first marginal $\tgX$ and $\Pi(*,\tgY)$ the set of couplings having second marginal $\tgY$, and $\Pi(\tgX,\tgY)=\Pi(\tgX,*)\cap\Pi(*,\tgY)$ the couplings with marginals $(\tgX,\tgY)$. We now recall an instrumental disintegration formula:%

Let $\pi,\bar\pi\in\cP(\X\times \Y)$ with $\pi \ll \bar\pi$, $ K_{\bar\pi}(x,dy)=\nicefrac{\bar\pi(dx,dy)}{\px \bar\pi(dx)}$. We have   $\bar\pi=\px \bar\pi \otimes K_{\bar\pi}$ and\footnote{The last equality uses: $\KL(
 \pi|\px  \pi \otimes K_{\bar\pi})=\int \ln\left(\frac{\px\pi\otimes K_{\pi}}{\px\pi \otimes K_{\bar{\pi}}}\right)d\px\pi\otimes K_{\pi} =\int_{\X} \KL(%
 K_{\pi}
 |K_{\bar\pi})\,d\px \pi%
 $
 .} 
\begin{equation}\label{eq:disintegration}
   \KL(\pi | \bar\pi)=\KL(\px \pi|\px \bar\pi)+\int_\X \KL(%
 K_{\pi}
 |K_{\bar\pi})\,d\px \pi =\KL(\px \pi|\px \bar\pi)+\KL(
 \pi|\px  \pi \otimes K_{\bar\pi}).
\end{equation}%
This decomposition is at the heart of the two objective functions $\Fsink$ and $\FEM$ considered below.

\subsection{Sinkhorn}

To describe the entropic optimal transport problem we follow~\cite{Nutz2021IntroductionTE}. Consider a cost function $c
\in L^\infty(\X\times \Y,\tgX\otimes\tgY)$
and a regularization parameter $\epsilon>0$. The entropic optimal transport problem is the minimization problem 
\begin{equation}\label{eq:ot_min_kl}
    \oteps(\tgX,\tgY)=\min_{\pi \in \Pi(\tgX,\tgY)}\KL(\pi|e^{-c/\eps}\tgX\otimes\tgY).
\end{equation}
By adding a constant to $c$ we can assume without loss of generality that $e^{-c/\eps}\tgX\otimes\tgY$ has mass $1$. Since $c$ is bounded, \eqref{eq:ot_min_kl} admits a unique solution $\pi_*$.
We use a characterization~\citep[Theorem 4.2, Lemma 4.9, Section 6]{Nutz2021IntroductionTE} of the set of cyclically invariant couplings, and define it as follows, as the set of couplings $\pi$ that solve an entropic optimal transport problem for their own marginals,
\begin{equation}
    \Pi_c =\{\pi \in \cP(\X\times \Y)\, | \, \KL(\pi|e^{-c/\eps}\mu\otimes\nu) = \min_{\pit\in\Pi(\mu,\nu)}\KL(\pit|e^{-c/\eps}\mu\otimes\nu),\, (\mu,\nu)=(\px \pi, \py \pi)\}.
\end{equation}
Moreover when $\pi \in \Pi_c$, there exist $f\in L^\infty(\X)$ and $g\in L^\infty(\Y)$ such that $\pi=e^{(f+g-c)/\eps}\mu\otimes\nu$.

The Sinkhorn algorithm in its primal formulation solves \eqref{eq:ot_min_kl} by alternating (entropic) projections on $\Pi(\tgX,*)$ and $\Pi(*,\tgY)$~\citep{Ruschendorf1995}, i.e.\ initializing with $\pi_0\in \PiexpInf$, iterate

\vspace{-8mm}

\begin{align}\label{eq:sinkhornprimal}
    \pi_{n+\frac12}&=\argmin_{\pi\in\Pi(\tgX,*)}\KL(\pi|\pi_n),\\
    \pi_{n+1}&=\argmin_{\pi\in\Pi(*,\tgY)}\KL(\pi|\pi_{n+\frac12}).\label{eq:sinkhornprimal-2}
\end{align}
Let $\mu_n=\px \pi_n$. More explicitly, \eqref{eq:sinkhornprimal} is a ``rescaling of the rows'', $\pi_{n+\frac12}(dx,dy)=\pi_n(dx,dy)\tg(dx)/\px\pi_n(dx)$, while~\eqref{eq:sinkhornprimal-2} is a ``rescaling of the columns'', $\pi_{n+1}(dx,dy)=\pi_{n+\frac12}(dx,dy)\tgY(dy)/\py\pi_{n+\frac12}(dy)$. This can be seen as a consequence of  \eqref{eq:disintegration}, as in \eqref{eq:sinkhornprimal} the first marginal is fixed, so the optimum of \eqref{eq:sinkhornprimal} is such that the integral term of \eqref{eq:disintegration} vanishes. One can also show recursively that $\pi_n\in\PiexpInf$ \citep[see][Section 6, Lemma 6.22]{Nutz2021IntroductionTE}.

Define the constraint set $C=\Pi(*,\tgY)$, $\cR=C\cap\PiexpInf$ and the objective function
\begin{equation}
    \Fsink(\pi)=\KL(\px \pi|\tgX).
\end{equation}

Connections between mirror descent and Sinkhorn iterations for the entropic regularized optimal transport problem were first discovered in \cite{mishchenko2019sinkhorn,menschpeyre, Leger2020}. We propose yet another mirror descent interpretation of Sinkhorn in the spirit of \citet[]{Leger2020}, and use the primal formulation \eqref{eq:sinkhornprimal}--\eqref{eq:sinkhornprimal-2} directly instead of introducing dual potentials. This is stated in the following Proposition, whose complete proof can be found in \Cref{sec:proof_prop_sinkhorn_mirror_descent}.

\begin{proposition}[Sinkhorn as mirror descent]\label{prop:sinkhorn_mirror_descent}
    The Sinkhorn iterations~\eqref{eq:sinkhornprimal} can be written as a mirror descent with objective $\Fsink$ and Bregman divergence $\KL$ over the constraint $C=\Pi(*,\tgY)$,
    \[
        \pi_{n+1}=\argmin_{\pi\in C} \ps{\nabla_{\!C} \Fsink(\pi_n),\pi-\pi_n} + \KL(\pi|\pi_n) \text{ with $\nabla_{\!C}\Fsink(\pi_n)= \ln(d\mu_n/d\tgX) \in L^\infty(\X\times\Y)$.}
    \]    
\end{proposition}
\begin{sproof} 
    Let $\mu_n=\px\pi_n$ where $\pi_n$ is defined in \eqref{eq:sinkhornprimal}. We have the identity:
    \begin{equation*}
        \Fsink(\pi_n)+\ps{\nabla_{\!C} \Fsink(\pi_n),\pi-\pi_n}+\KL(\pi|\pi_{n})
        =\KL(\pi|\tgX\otimes\nicefrac{\pi_n}{\mu_n})=\KL(\pi|\pi_{n+\frac12}).
    \end{equation*}
    We conclude by taking the argmin over $\pi\in C$.
\end{sproof}
We first show relative smoothness of $\Fsink$ relatively to $\phine$, as a consequence of the standard KL data processing inequality, i.e.\ KL of the marginals is smaller than Kl of the plans.
\begin{lemma}\label{prop:rel-smooth-Sink}
    The functional $\Fsink$ is convex and is $1$-relatively smooth w.r.t. $\phine$ over $\cP(\X\times\Y)$. 
\end{lemma}
\begin{proof}
    Let $\pi,\pit\in\cP(\X\times\Y)$ with $\px \pit\ll \px \pi \ll \tgX$. Then with straightforward computations, $D_{\Fsink}(\pit|\pi)=\KL(\px \pit|\px \pi)\ge 0$, so $\Fsink$ is convex. Then  \eqref{eq:disintegration} results in  $D_{\Fsink}(\pit|\pi)\le\KL(\pit|\pi)$.%
\end{proof}
By considering singular first marginals, it is obvious that there exists no $l>0$ such that $D_{\Fsink}(\pit|\pi)\ge l\KL(\pit|\pi)$ for all $\pit,\pi\in C$; i.e.\ that relative strong convexity of $\Fsink$ relatively to $\phine$ does not hold \textit{over all $C$}. However, we show in the next Proposition that this inequality actually holds over $\cR=\PiexpInf\cap C$. 
Its complete proof can be found in \Cref{sec:proof_quantitative-stability}.

\begin{proposition} \label{prop:quantitative-stability}%
    Let $D_c:=\frac12 \sup_{x,y,x',y'}[c(x,y)+c(x',y')-c(x,y')-c(x',y)]<\infty$. For $\pit,\pi\in \PiexpInf\cap C$, we have that
    \begin{equation}\label{eq:Fsink_strong_conv}
        \KL(\pit|\pi)\le (1+4 e^{3D_c/\epsilon}) \KL(\px \pit |\px \pi),
    \end{equation}
    in other words $\Fsink$ is $(1+4 e^{3D_c/\epsilon})^{-1}$-relatively strongly convex w.r.t. $\KL$ over $\PiexpInf\cap C$.
\end{proposition}
\begin{sproof}
    For $\pi,\pit\in\PiexpInf\cap C$ with their potentials and marginals $(f,g,\mu,\tgY)$ and $(\tilde f,\tilde g,\mut,\tgY)$ respectively, setting $\normvar{f}=(\sup_\X f)-(\inf_\X f)$, we can derive the bound
    \begin{equation}\label{eq:kl_bound}
    \eps\KL(\pit|\pi)\le \normvar{\tilde f-f}\normtv{\mut-\mu} +\eps\KL(\mut|\mu).
    \end{equation}
  We can then bound the potentials  by the marginals using \Cref{quantitative-stability-estimate} in Appendix:
    \begin{equation}\label{eq:potentials_bound}
    \normvar{f-\ft}+\normvar{g-\gt}\le 2 \eps\,e^{3D_c/\epsilon}\big(\normtv{\mu-\mut} + \normtv{\tgY-\tgY}\big).
        \end{equation}
    Then, chaining~\eqref{eq:kl_bound} and~\eqref{eq:potentials_bound}  %
    we conclude using Pinsker's inequality.
\end{sproof}
 We are now ready to recover convergence rates for Sinkhorn leveraging relative smoothness and strong convexity.

\begin{proposition}[Sinkhorn convergence] %
For all $n\ge 0$, the Sinkhorn iterates verify, for $\pi_*$ the optimum of \eqref{eq:ot_min_kl} and $\mu_*$ its first marginal,
    \begin{equation}\label{eq:rate_sinkhorn}
        \KL(\mu_n|\mu_*) \le\frac{  \KL(\pi_*|\pi_0)}{(1+4e^{\frac{3Dc}{\epsilon}})\left(\left( 1+ 4e^{-\frac{3 D_c}{\epsilon} } \right)^n -1\right) }\le \frac{\KL(\pi_*|\pi_0)}{ n}.
    \end{equation}
\end{proposition}
\begin{proof} Fix $n\ge 0$, we know that $\pi_n,\pi_*\in \cR:=C\cap \PiexpInf$. We conclude by applying \Cref{th:rate} to $\cF=\Fsink$ and $\phi=\KL(\cdot | \pi_*)$, leveraging the results of \Cref{prop:rel-smooth-Sink} and \Cref{prop:quantitative-stability}.
\end{proof}
The linear convergence of Sinkhorn for bounded costs $c$ has been known since at least~\citet[]{Franklin1989} and has then been derived also in the non-discrete case and in multimarginal settings \citep[see][and references therein]{Carlier2022OnTL}. To derive this linear rate (the first inequality in \eqref{eq:rate_sinkhorn}), we require relative strong convexity of the objective - hence we also fundamentally rely on the boundedness of the cost, as one can see from the assumptions of \Cref{prop:quantitative-stability}. Indeed the proof of \Cref{prop:quantitative-stability} relies on the classical result that the soft $c$-transforms are contractions in the Hilbert metric; a result also at the heart of other proofs for the linear rate~\citep[see][for a proof]{Franklin1989,chen2016entropic}.
Regarding the sublinear convergence (the second inequality in \eqref{eq:rate_sinkhorn}), \citet{Leger2020} first obtained sublinear rates for unbounded costs leveraging relative smoothness, using \eqref{eq:dual_iteration} formally and through dual iterations on the potentials. In this Section we assumed $c\in L^{\infty}(\X \times \Y)$, in order to manipulate finite quantities in our computations (e.g. first variations in the proof of \Cref{prop:sinkhorn_mirror_descent}) - hence the latter can be seen as a convenient working hypothesis. 
In this paper we derive the same rate as \citet{Leger2020} rigorously with a more direct proof using primal iterations, and complete the picture by recovering linear rates of convergence.

\subsection{ Expectation-Maximization}

In this subsection, we show how the EM algorithm can always formally be written as a mirror descent scheme, and, when optimizing the latent variable distribution, results in a convex problem with sublinear convergence rates. Consider the following probabilistic model: we have a latent, hidden random variable $X\in (\X,\tgX)$, an observed variable $Y\in \Y$ distributed as $\tgY$, and we posit a joint distribution $p_{q}(dx,dy)$ parametrized by an element $\qEM$ of some given set $\Q$. As presented in \citet[]{Neal1998}, the goal is to infer $q$ by solving
\begin{equation}
    \min_{\qEM\in\Q} \KL(\tgY|\py p_\qEM),
\end{equation}
where $\py p_{\qEM}(dy)=\int_{\X} p_{\qEM}(dx,dy)$. The EM approach starts by minimizing a surrogate function of $\qEM$ upperbounding $\KL(\tgY|\py p_\qEM)$. For any $\pi\in\Pi(*,\tgY)$, by the data processing inequality,
\[
    \KL(\tgY|\py p_\qEM)\le \KL(\pi|p_\qEM)=:L(\pi,\qEM).
\]
Again, as a consequence of the disintegration formula \eqref{eq:disintegration}, there is  equality if and only if
\begin{equation}\label{eq:EM_iff}
    \pi(dx,dy)=p_\qEM(dx,dy)\tgY(dy)/\py p_\qEM(dy).
\end{equation}

EM then proceeds by alternate minimizations of $L(\pi,\qEM)$ \citep[see][Theorem 1]{Neal1998}:%
\begin{align} 
\label{eq:em-iterations-M-step}
    \qEM_{n} &= \argmin_{\qEM \in \Q} \KL(\pi_{n}|p_{\qEM}),\\
    \label{eq:em-iterations-E-step}
    \pi_{n+1}&=\argmin_{\pi\in\Pi(*,\tgY)}\KL(\pi|p_{\qEM_n}).
\end{align}
\vspace{-5mm}

The above formulation consists in \eqref{eq:em-iterations-M-step}, optimizing  the parameters $q_n$ at step $n$ (M-step), and then \eqref{eq:em-iterations-E-step}, optimizing  the joint distribution $\pi_{n+1}$ at step $n+1$ (E-step). We choose this order to highlight the analogy of  \eqref{eq:em-iterations-E-step} with Sinkhorn's \eqref{eq:sinkhornprimal-2}. The minimization \eqref{eq:em-iterations-E-step} corresponds to taking an explicit expectation, according to \eqref{eq:EM_iff} justifying the denomination. On the contrary making explicit the M-step is often difficult.

Define the constraint set $C=\Pi(*,\tgY)$ and the, possibly non-convex, objective function
\begin{equation}\label{eq:def-FEM}
  \FEM(\pi) = \inf_{\qEM \in \Q}\KL(\pi|p_\qEM).  
\end{equation}

We first show that EM can be formally written as a mirror descent scheme in the following Proposition, whose precise statement with additional assumptions can be found in \Cref{sec:proof_enveloppe_FEM}. %

\begin{proposition}[EM as mirror descent, formal]\label{prop:em-mirror-descent}The EM iterations~\eqref{eq:em-iterations-M-step}--\eqref{eq:em-iterations-E-step} can be written as a mirror descent iteration with objective function $\FEM$, Bregman potential $\phine$ and constraints $C$,

\vspace{-3mm}
    \begin{equation}\label{eq:em_md_iterate}
        \pi_{n+1}=\argmin_{\pi\in C} \ps{\nabla_{\!C} \FEM(\pi_n),\pi-\pi_n} + \KL(\pi|\pi_n) \text{ with $\nabla_{\!C}\FEM(\pi_n)= \ln(d\pi_n/dp_{\qEM_n}).$}
        \end{equation}
\end{proposition}
\vspace{-3mm}
\begin{sproof}
 Let $\pi_n$ be the current EM iterate. Formally, we use the envelope theorem to differentiate $\FEM$ and find that $\nabla_{\!C}\FEM(\pi_n)=\ln(d\pi_n/dp_{\qEM_n})$ (see \Cref{sec:proof_enveloppe_FEM} for a justification based on directional derivatives and \citet[Theorem 3]{Milgrom02envelopetheorems}). Then for any coupling $\pi$, we have the identity
\begin{align*}
  &\FEM(\pi_n)+\ps{\nabla_{\!C}\FEM(\pi_n),\pi-\pi_n}+\KL(\pi|\pi_n) %
 = \KL(\pi|p_{\qEM_n}).
\end{align*}
Thus \eqref{eq:em_md_iterate} matches~\eqref{eq:em-iterations-E-step}.
\end{sproof}
Since $\FEM$ is in general non-convex, we cannot apply outright the framework developed in \Cref{sec:MD_convergence}. There is however one direct case of $p_q$ making $\FEM$ convex, by optimizing only over its first marginal.

\textbf{Latent EM.} We consider the case where $p_q(dx,dy)=\mu(dx)K(x,dy)$, i.e. $p_q$ is of the form $\mu\otimes K$, with $\Q=\cP(\X)$ and $K$ kept fixed along iterations. In other words, we choose to only optimize over the density of the latent variable, and keep the mixture parameters fixed. Here $K$ can be interpreted as the conditional distribution of $Y$ given $X$, $K(x,dy)=\mathbb{P}(Y=y\mid X=x)$. We consider in the following Gibbs distributions with $K(x,dy)=e^{-c(x,y)}\tgY(dy)$ with $c$ uniformly bounded. The term $c$ can be interpreted as a cost similarly to the entropic optimal transport \eqref{eq:ot_min_kl}. 

\begin{remark}[Various EM]%
The general goal of EM is to fit, through the objective function $F_{\text{EM}}$, a parametric distribution, e.g. a mixture of Gaussians,  to some observed data $Y$. One needs to estimate both the latent variable distribution on $X$  (i.e. weights of each Gaussian) and the parameters of conditionals $P(Y|X=x)$ (e.g. means and covariances of each Gaussian). Latent EM focuses on learning the mixture weights, since it consists in optimizing over the nonparametric latent distribution $\mu$, which can be continuous or discrete. In contrast, the parametric setting considered by \citet[]{Kunstner2021Homeomorphic}, who obtained $\mathcal{O}(1/n)$ rates of convergence in KL for EM, can be seen as complementary to ours since they consider a fixed $\mu$ and a variable exponential mixture $K_\theta$, with $q=\theta$. 
\end{remark}

Parametrizing by the first marginal, EM iterations  \eqref{eq:em-iterations-M-step}--\eqref{eq:em-iterations-E-step} takes the following form for Latent EM:
\begin{align}
\label{eq:em-latent-iterations-M-step}
 \mu_{n} &= \argmin_{\mu \in \cP(\X)} \KL(\pi_{n}|\mu\otimes K),\\
 \label{eq:em-latent-iterations-E-step}
    \pi_{n+1}&=\argmin_{\pi\in\Pi(*,\tgY)}\KL(\pi|\mu_n\otimes K).
\end{align}

First, we necessarily have from \eqref{eq:em-latent-iterations-M-step} that $\mu_{n} = \px \pi_{n}$. Indeed, from the disintegration formula \eqref{eq:disintegration}, \eqref{eq:em-latent-iterations-M-step} corresponds to minimizing over first marginals. Then, since the E-step \eqref{eq:em-latent-iterations-E-step} corresponds to computing \eqref{eq:EM_iff}, we can rewrite \eqref{eq:em-latent-iterations-M-step}--\eqref{eq:em-latent-iterations-E-step} as:

\vspace{-4mm}

    \begin{equation}\label{eq:em-latent}
    \mu_{n+1}(\cdot) = \int_{\Y} \pi_{n+1}(\cdot,dy) = \mu_{n}(\cdot)
\int_{\mathcal{Y}} \frac{K(\cdot,dy)\tgY(dy)}{\int_{\mathcal{X}} K(x,dy)\mu_{n}(dx)}.
\end{equation}
\vspace{-2mm}

Define $\FEMK(\pi):=\inf_{\mu\in \cP(\X)}\KL(\pi|\mu\otimes K)$. Notice that by disintegration \eqref{eq:disintegration},  $\FEMK$ takes the form $\FEMK(\pi)=\KL(\pi|\px\pi\otimes K)$. 
To take care of the initialization, we define the operator $T_K: \mu\in \cP(\X)\mapsto \int_{\X} \mu(dx)K(x,\cdot) \in \cM_r(\Y)$, with $K(x,dy)=k(x,y)\tgY(dy)$ and take $\mu_0=e^{f_0}\tgX$ with $f_0\in L^\infty$ and assume that $T_K \tgX \gg \tgY$ (in other words we assume that the mixture applied to the latent space is compatible with all the observations).

\begin{proposition}[Latent EM as mirror descent]\label{prop:em_mirror_descent}

   The latent EM iterations \eqref{eq:em-latent-iterations-M-step}--\eqref{eq:em-latent-iterations-E-step} can be written as mirror descent with objective $\FEMK$, Bregman potential $\phine$ and the constraints $C=\Pi(*,\tgY)$,
    \begin{equation*}
        \pi_{n+1}=\argmin_{\pi\in C} \ps{\nabla_{\!C} \FEMK(\pi_n),\pi-\pi_n} + \KL(\pi|\pi_n) \text{ with $\nabla_{\!C}\FEMK(\pi_n)= \ln\left(\frac{d\pi_n}{d(\mu_n\otimes K)}\right) \in L^\infty$.}
    \end{equation*}
\end{proposition}

\vspace{-5mm}
\begin{proof} Similarly to \Cref{prop:em-mirror-descent}, 
we have the identity
    \[
        \FEMK(\pi_n)+\ps{\nabla_{\!C} \FEMK(\pi_n),\pi-\pi_n}+\KL(\pi|\pi_{n}) = \KL(\pi|\px \pi_n\otimes K),
    \]
    where $\nabla_{\!C}\FEMK(\pi_n)=\ln(\pi_n/\px \pi_n \otimes K) \in L^\infty(\X,\R)$, see \Cref{sec:proof_first_variations_FEMK} for rigorous justifications. Since $\mu_n=\px \pi_n$ due to \eqref{eq:em-latent-iterations-M-step}, we conclude by minimizing over $\pi\in C$. 
\end{proof}
We are now ready to state convergence rates of latent EM in the following proposition. The reader may refer to \Cref{sec:proof_rates_em} for a complete proof.
\begin{proposition}[Convergence rate for Latent EM]\label{prop:rates_em}
Set $\mu_*\in \argmin_{\mu\in\cP(\X)} \KL(\tgY | T_K(\mu))$. The functional $\FEMK$ is convex and $1$-smooth relative to $\phine$. Moreover for $\pi_0\in \Pi(*,\tgY)$,
    \[
        \KL(\tgY|T_K\mu_n)\le \KL(\tgY|T_K\mu_*)+ \frac{\KL(\mu_*|\mu_0)+\KL(\tgY|T_K\mu_*)-\KL(\tgY|T_K\mu_0)}{n}.
    \]
\end{proposition}
\vspace{-5mm}

\begin{sproof}
By the disintegration formula \eqref{eq:disintegration}, we can decompose $\FEMK$ as
\begin{equation}\label{eq:decompose_FEMK}
    \FEMK(\pi) = \KL(\tgY|\py (\px \pi \otimes K)) +\int \KL(\pi/\tgY| (\px \pi \otimes K )/ \py (\px \pi \otimes K) ) d\tgY,
\end{equation}
 and show then straightforwardly that $\pi_*$ defined by $\pi_*(dx,dy)=\mu_*(dx)k(x,dy)\tgY(dy)/(T_K\mu_*)(dy)$ is a minimizer of $\FEMK$ with $\FEMK(\pi^*) = \KL(\tgY|T_K \mu_*)$. %
 Then, by the disintegration formula \eqref{eq:disintegration} and linearity of the Bregman divergence, $\KL(\pi|\tilde{\pi})= D_{\Fsink}(\pi|\tilde{\pi}) + D_{\FEMK}(\pi|\tilde{\pi})$, hence $\FEMK$ is 1-relatively smooth w.r.t. $\phine$. Consequently,  \Cref{th:rate} yields:
  \[
        \FEMK(\pi_n)\le \FEMK(\pi_*)+\frac{\KL(\pi_*|\pi_0)}{n},
    \]
 Finally: $\KL(\tgY|T_K \mu_n)=\KL(\py\pi_n|\py(\px\pi_n \otimes K))\le \KL(\pi_n|\px\pi_n \otimes K)= \FEMK(\pi_n)$.
\end{sproof}

\begin{remark}[Richardson--Lucy] Interestingly, the iterations \eqref{eq:em-latent} of latent EM correspond precisely to that of Richardson--Lucy deconvolution \citep[][]{Richardson1972,Lucy1974} where $K$ is a known convolution and one aims at recovering the original signal $\mu_*$ based on the observations $\tgY$. Thus our proof yields 
rates of convergence for this other algorithm from signal processing, a novel result to the best of our knowledge.
\end{remark}

\tb{Conclusion:} We have provided a rigorous proof of convergence of mirror descent under relative smoothness and convexity, which holds in the infinite-dimensional setting of optimization over measure spaces. The latter condition can handle objective functionals that are not smooth in the standard sense, such as the ubiquitous KL. It enabled us to provide a new and simple way to derive rates of convergence for Sinkhorn's algorithm. We also derived new convergence rates for EM when restricted to the latent distribution, obtaining complementary rates to \cite{Kunstner2021Homeomorphic}.

\newpage
\bibliographystyle{apalike}
\bibliography{biblio}

\appendix
\onecolumn

\section*{Appendix}

\section{Definition of Gâteaux and Fréchet derivatives}\label{sec:add_defs}
We first recall the notion of Gâteaux and Fréchet derivatives for $\cF : \cM(\X ) \rightarrow \R \cup \left\{\pm\infty\right\}$ where $ \cM(\X )$ is a topological vector space \citep[Chapter 7, pp.267,273]{aliprantis2006infinite}, see also \citet[Section 1]{Phelps1989} for Banach spaces.
\begin{definition} The function
$\cF$ is said to be Gâteaux differentiable at $\nu$ if there exists a linear operator $\nabla F(\nu):\cM(\X ) \rightarrow \R$ such that for any direction $\mu\in \cM(\X )$:
\begin{equation}\label{eq:gateaux}
   \nabla \cF(\nu)(\mu)= \lim_{h\rightarrow 0}\frac{\cF(\nu+h\mu)-\cF(\nu)}{h}.
\end{equation}
The operator $\nabla \cF(\nu)$ is called the Gâteaux derivative of $\cF$ at $\nu$, and if it exists, it is unique.
\end{definition}
\begin{definition} If $\cM(\X)$ is a normed space, the function
$\cF$ is said to be Fréchet differentiable at $\nu$ if there exists a bounded linear form $\delta\cF(\nu,\cdot):\cM(\X ) \rightarrow \R$ such that
\begin{equation}
    \cF(\nu+h\mu)=\cF(\mu)+h\delta \cF(\nu,\mu)+h o(\|\mu\|_{\cM(\X)})
\end{equation}
Equivalently, the operator $\delta \cF(\nu,\cdot)$ is called the Fréchet derivative of $\cF$ at $\nu$ if it is a Gâteaux derivative of $\cF$ at $\nu$ and the limit \eqref{eq:gateaux} holds uniformly
in $\mu$ in the unit ball (or unit sphere) in $\cM(\X)$.
\end{definition}
If $\cF$ is Fréchet differentiable, then it is also Gâteaux differentiable, and its Fréchet and Gâteaux derivatives agree: $\nabla \cF(\nu)(\mu)=\delta \cF(\nu,\mu)$.

\section{Additional details on the well-posedness of the mirror descent scheme}\label{sec:assumptions_details}

Recall \Cref{ass:diff+inf-compact}(Lower semicontinuity and coercivity):  (i) the set $C$ is $\tau$-closed in $\cM(\X)$, the functionals $\cG_n(\cdot):= d^+\!\cF(\mu_n)(\cdot- \mu_n)$ and $D_{\phi}(\cdot | \mu_n)$ are proper and $\tau$-l.s.c.\ when restricted to $C$, and the functional $\cG_n+D_{\phi}(\cdot |\mu_n)+i_{C}$\footnote{$i_{C}$ denotes the indicator function of the set $C$, defined by $i_{C}(\mu)=0$ if $\mu\in C$, $+\infty$ otherwise for any $\mu\in \cM(\X)$. Notice that $i_{C}$ being $\tau$-l.s.c.\ is equivalent to $C$ being $\tau$-closed in $\cM(\X).$} has at least one $\tau$-compact sublevel set. (ii) For each $n\ge0$, the first variations $\nabla_{\!C}\phi(\mu_n)$ exist. (iii) The iterates belong to $\cR$.

\textbf{Case where $\cF$ has first variations.} Equip $\cM(\X)$ with a topology $\tau$ that is stronger than the $\cM(\X)^*$-weak topology. If $\cF$ has first variations, then we can even remove $\cG_n(\cdot)$ from \Cref{ass:diff+inf-compact}, since $\cG_n(\cdot)$ is linear on the set of interest and $\tau$-l.s.c. Whence we get the simpler assumption
\begin{assumplist3}
	\setlength\itemsep{0.2em}
	\item \label{ass:diff+inf-compact_simple}(Lower semicontinuity and coercivity) For each $n\ge0$, the iterates belong to $\cR$ and the first variations $\nabla_{\!C}\phi(\mu_n)$ and $\nabla_{\!C}\cF(\mu_n)$ exist. Moreover the set $C$ is $\tau$-closed in $\cM(\X)$, the functional $\phi$ is proper and $\tau$-l.s.c.\ when restricted to $C$, and $\phi(\cdot)$ has at least one $\tau$-compact sublevel set when restricted to $C\cap\dom(\cF)$.%
\end{assumplist3}

\textbf{Weakly compact sets of $\cM(\X)$.}
In finite dimensions, a set is compact iff bounded and closed; however, in infinite dimensions, characterizing compact sets is more delicate. Below we recall some classical set of conditions that guarantee (weak) compactness or lower semicontinuity.

If $\cM(\X)$ is a reflexive Banach space, then the weakly compact sets are just the bounded weakly closed sets, as a consequence of the Banach--Alaoglu theorem \citep[see e.g.][Theorem 2.4.2]{attouch2014variational}. In other cases, one needs more specific theorems such as Dunford--Pettis' theorem for $L^1(\rho)$ \citep[see e.g.][Theorem 2.4.5]{attouch2014variational}. Since we are dealing with convex functions, for normed $\cM(\X)$, strongly closed sublevel sets are also weakly closed, a result known as Mazur's lemma. So the notions of weakly l.s.c.\ and strongly l.s.c.\ convex functions coincide, as recalled in \citet[Theorem 3.3.3]{attouch2014variational}.

We now regroup some known properties of KL, in particular to show that \Cref{ass:diff+inf-compact} holds for $\phi=\phine$.
\begin{remark}[Properties of KL]\label{rmk:KL_properties_app}
 For compact $\X$, the domain of the negative entropy $\phine$ is strictly included in $L^1_+(\X)$, contains $L^q_+(\X)$ for $q>1$, and is of empty interior for the norm/strong topology of $L^q(\X)$ \cite[Lemma 4.1]{resmerita2005regularization}. Regarding the use of $\phi=\phine$ in \eqref{eq:Bregman_prox_algo}, one can for instance take $\cM(\X)= L^1(\X)$ equipped with the weak topology induced by $L^\infty(\X)$ and the Lebesgue measure as reference. We have that $\cP(X)\cap L^1(\X)$ is weakly closed and that $\KL$ and $\phine$ are strictly convex, weakly l.s.c. and have weakly compact sublevel sets in $L^1(\X)$ by \citet[Lemma 2.1, 2.3]{eggermont1993maximum} \citep[see also][Section 3]{Resmerita2007}. By \citet[Lemma 4.1]{resmerita2005regularization}, a sufficient condition for $\KL$ (resp.\ $\phine$) to have a first variation in $L^\infty$ at $\mu$ is that there exists $\kappa_0,\kappa_1>0$ such that $\kappa_0\le\nicefrac{d\mu}{d\tg}(x)\le\kappa_1$ almost everywhere over $\X$ (resp. $\phine$ for $\tg=\rho$). $\KL$ is not Gâteaux-differentiable for non-finite $\X$ as recalled for instance in \citet[p12]{butnariu2006bregman} and \citet[Remark 7.13][p239]{Santambrogio2015}.
\end{remark}

As a follow-up of \Cref{rmk:proof_cv_main}, we now give some known conditions for a sum of subdifferentials to be the subdifferential of the sum.
\begin{remark}[About the proof of convergence in \Cref{th:rate}]\label{rmk:proof_cv_app}
Our proof of \Cref{th:rate} resembles the one of \citet{lu2018relatively}, which also relies on a three-point inequality as stated in  \Cref{lem:three-point}. However, the proof of the latter inequality in finite dimensions uses a formula of the form $\partial (\cG+D_\phi) = \partial \cG +\partial D_\phi$ as in \citet[Lemma 3.2]{chen1993convergence}, but which is harder to derive in infinite dimensions.
Such an equality between subdifferentials can be obtained typically under at least three (non-equivalent) conditions for convex and l.s.c.\ $\cG$ and $\phi$ over a Banach space $\cM(\X)$: (i) having $\cup_{\lambda \ge 0}\lambda (\dom(\cG) -\dom(D_\phi))$ to be a closed vector space of $\cM(\X)$ \citep{attouch1986duality}; (ii) having a non-empty (quasi) relative interior of $(\dom(\cG)-\dom(D_\phi))$ \citep{borwein2003notions}%
; (iii) continuity of $D_\phi$ or $\cG$ at least at some $\mu\in \dom(\cG)\cap\dom(D_\phi)$ \citep[Theorem 3.30]{peypouquet2015convex}. Condition (iii) does not hold when $\cF$ and $D_{\phi}$ are chosen as the KL divergence (defined below in \Cref{rmk:KL_properties}) in none of the spaces we consider since $\KL$ is not continuous, its domain being of empty interior. The other conditions are difficult to verify for given functionals. For instance, $\dom(\phi)$ is not explicit for the negative entropy (see \Cref{rmk:KL_properties}).
On the contrary, by favoring directional derivatives and first variations, we circumvent  most of the difficulties related to (sub)differentiability.
\end{remark}

\section{Additional technical results}

	\begin{lemma}\label{lem:conv_directionnal_deriv} Let $f$ be a proper function over a vector space $\Y$ with values in $\R\cup\{+\infty\}$. The following conditions are equivalent:
	\begin{enumerate}[label=\roman*)]
		\item $f$ is convex;
		\item $\dom(f)$ is convex, and, for all $x,y\in \dom(f)$, $d^+\!f(x)(y-x)$ exists, with value in $\R\cup\{-\infty\}$, and we have $f(x)+	d^+\!f(x)(y-x) \le f(y)$, i.e.\ $D_f(y|x)\ge 0$;
		\item $\dom(f)$ is convex, and, for all $x,y\in \dom(f)$, $d^+\!f(x)(y-x)$ exists, with value in $\R\cup\{-\infty\}$, and we have $d^+\!f(x)(y-x)+d^+\!f(y)(x-y)\le 0$.
	\end{enumerate}
	\end{lemma}
	The lemma immediately extends to strictly convex functions by taking strict inequalities.
	\begin{proof}
	    Given $x,y\in \dom(f)$, define for any $\lambda\in[0,1]$ $u_\lambda=x+\lambda(y-x)$. Assuming (i), then
	    \begin{gather*}
	        f(u_\lambda)\le \lambda f(y)+(1-\lambda)f(x)\\
	        f(x)+\frac{f(u_\lambda)-f(x)}{\lambda} \le f(y),
	    \end{gather*}
	    which yields (ii) by the decreasingness of differential quotients discussed in \Cref{rmk:directional_derivative}.

	    Assuming (ii), we just sum the two inequalities $(f(x)+	d^+\!f(x)(y-x) \le f(y))$ and $(f(x)+	d^+\!f(x)(y-x) \le f(y))$, to derive (iii).
	    
	    The last implication to show is (iii)$\Rightarrow$ (i) which requires to perform an integration. Assume that (iii) holds, we want to show that $f(u_\lambda)\le \lambda f(y)+(1-\lambda)f(x)$. Set $g(\lambda):=f(u_\lambda)$ and denote by $g_+'(\lambda)$ (resp.\ $g_-'(\lambda)$) its right (resp.\ left) derivative, both derivatives exist with value in $\R\cup\{-\infty\}$ for any $\lambda\in(0,1)$ since
	    \begin{align*}
	        d^+\!f(u_\lambda)(y-x)=\lim_{h \rightarrow 0^+}\frac{f(u_\lambda+h(y-x))- f(x)}{h}=g_+'(\lambda)
	    \end{align*}
	    similarly $d^+\!f(u_\lambda)(x-y)=-g_-'(\lambda)$. Consequently, for all $0<\lambda_1<\lambda_2<1$, applying (iii) to $u_{\lambda_1}$ and $u_{\lambda_2}$, we have that $g_+'(\lambda_1)\le g_-'(\lambda_2)$. We now show that $\lambda\mapsto g_-'(\lambda)$ is increasing over $(0,1)$. We just have to show that $g_-'(\lambda_1)\le \sup_{\lambda\in(0,\lambda_2)}g_+'(\lambda)$. By contradiction, we could fix $\lambda_1\in(0,\lambda_2)$ and $\epsilon>0$ such that, for all $\lambda\in(0,\lambda_2)$,  $g_-'(\lambda_1)\ge g_+'(\lambda) +\epsilon$. By definition of the directional derivatives, we can then fix $\delta_1\in (0,\lambda_1)$ and $\lambda\in(\lambda_1-\delta_1,\lambda_1)$ such that for all $h_0\in(0,\delta_1)$
	    \begin{align*}
	        |g_-'(\lambda_1)+\frac{g(\lambda_1-h_0)- g(\lambda_1)}{h_0}|\le \epsilon/4 \quad ; \quad |g_+'(\lambda)-\frac{g(\lambda_1)- g(\lambda)}{\lambda_1-\lambda}|\le \epsilon/4 
	    \end{align*}
	     whence
	    \begin{align*}
	        \frac{g(\lambda_1)-g(\lambda_1-h_0)}{h_0}\ge \frac{g(\lambda_1)- g(\lambda)}{\lambda_1-\lambda}+\epsilon/2
	    \end{align*}%
	    which leads to a contradiction for $h_0=\lambda_1-\lambda$. Therefore $\lambda\rightarrow g_-'(\lambda)$ is increasing over $(0,1)$, upper bounded by $g_-'(y)=d^+\!f(u_1)(x-y)$. Since $g$ has both left and right derivatives, it is continuous over $[0,1]$. We can now apply (iii) to $x$ and $u_{\lambda}$, use the positive homogeneity of the directional derivative (which always holds by definition), and integrate over $(0,1)$ since the function $g_-'$ is Riemann-integrable,
	    \begin{align*}
	        0 & \ge d^+\!f(x)(u_{\lambda}-x)+d^+\!f(u_{\lambda})(x-u_{\lambda}) \\
	        & = \lambda d^+\!f(x)(y-x)+\lambda d^+\!f(u_{\lambda})(x-y) \\
	        0 & \ge d^+\!f(x)(y-x)-\int_0^1 g_-'(\lambda) d\lambda \\
	        & =d^+\!f(x)(y-x)+g(0)-g(1)  \\
	         & =d^+\!f(x)(y-x)+f(x) - f(y),
	    \end{align*}
	    which concludes the proof.
	\end{proof}
	
Below, we derive some useful characterizations of relative smoothness, 
by analogy with  \citet[][Proposition 1]{bauschke2017descent} for differentiable functions in finite dimensions. Similar results hold for relative convexity by the same arguments.

\begin{lemma}\label{lem:equivalences}
The following conditions are equivalent:
\begin{enumerate}[label=(\roman*),wide, labelwidth=0pt, labelindent=2pt,itemsep=-1ex,topsep=0pt]
		\item $\cF$ is $L$-smooth relative to $\phi$ over $C$;
		\item $L \phi- \cF$ is convex on $C\cap \dom(\phi) \cap \dom(\cF)$;
	\end{enumerate}
	and, if the first variations of $\cF$ and $\phi$ over $C$ evaluated at $\mu,\nu\in C\cap \dom(\phi) \cap \dom(\cF)$ exist,
	\begin{enumerate}[label=(\roman*),wide, labelwidth=0pt, labelindent=2pt,itemsep=-1ex,topsep=0pt]
		\item[(iii)]
		$\ps{\nabla_{\!C} \cF(\mu)-\nabla_{\!C} \cF(\nu), \mu - \nu}\le L \ps{\nabla_{\!C} \phi(\mu)-\nabla_{\!C} \phi(\nu),\mu-\nu}$
	\end{enumerate}
\end{lemma}
\begin{proof}
	This is a consequence of \Cref{lem:conv_directionnal_deriv} applied to $\psi(\mu)=L\phi(\mu)-\cF(\mu)$. More precisely, condition (i) can be written as
	$	d^+\psi(\mu)(\nu-\mu)\le \psi(\nu)-\psi(\mu) $
    which is equivalent to the convexity of $\psi$ by \Cref{lem:conv_directionnal_deriv}, hence (i)$\Leftrightarrow$(ii).
	Provided the first variations of $\cF$ and $\phi$ over $C$ exist, assuming (i) and (iii) boils down to \Cref{lem:conv_directionnal_deriv}-iii). Conversely, assuming (iii), we use \Cref{lem:conv_directionnal_deriv}-iii) and the linearity of the first variation \eqref{eq:first_var_dual}.
\end{proof}

\section{Smoothness of the Maximum Mean Discrepancy relatively to the KL divergence}\label{sec:app_mmd_kl}
Let $k : \X \times \X \to \R$ be a positive semi-definite kernel, $\kH$ its corresponding Reproducing Kernel Hilbert Space \citep{steinwart2008support}. The space $\kH$ is a Hilbert space with inner product and norm $\Vert \cdot \Vert_{\kH}$ satisfiying the reproducing property: for all 
$f \in \kH \text{ and }x\in \X,\; f(x)=\ps{f,k(x,\cdot)}_{\kH}$. For any $\mu \in \cP(\X)$ such that $\int \sqrt{k(x,x)}d\mu(x)<\infty$, the kernel mean embedding of $\mu$, $m_{\mu}= \int k(x,\cdot)d\mu(x)$, is well-defined, belongs to $\kH$, and $\E_{\mu}[f(X)] =\ps{f, m_{\mu}}_{\kH}$ \citep{smola2007hilbert}.
The kernel $k$ is said to be characteristic when such mean embedding is injective, that is, when any probability distribution is associated to a unique mean embedding. In this case, the kernel defines a distance between probability distributions referred to as the Maximum Mean Discrepancy (MMD), defined through the square norm of the difference between mean embeddings:
\begin{equation*}
	\MMD^2(\mu,\tg)  =\Vert m_{\mu}-m_{\tg} \Vert^2_{\kH} %
\end{equation*}
Interestingly, as soon as the kernel is bounded, the MMD is relatively smooth with respect to $\phine$, see \Cref{prop:rel_smooth_mmd} below. Notice that, thanks to the reproducing property,  $\mu \mapsto \MMD(\mu,\tg)$ is strictly convex whenever the kernel $k$ is characteristic, as it is the case for the Gaussian kernel. Similarly to $\KL$, the $\MMD$ can be written as a Bregman divergence of $\phi_k(\mu)=\|m_{\mu}\|_{\kH}^2=\int k(x,x')d\mu(x)d\mu(x')$.

\begin{proposition}\label{prop:rel_smooth_mmd} Let $\phine:\mu \mapsto \int \log(\mu)d\mu$ and fix $\nu \in \cP(\X)$. Take $k:\X\times\X\rightarrow\R$ to be a bounded semipositive definite kernel, i.e.\ $c_k = \sup_{x\in \X}k(x,x)<\infty$. The squared Maximum Mean Discrepancy $\MMD^2(\cdot,\nu)$ is $4c_k$-smooth relative to $\phine$.
\end{proposition}
\begin{proof} Let $\mu,\nu \in \cP(\X)$ and $f_{\mu,\tg}=\int k(x,\cdot)d\mu(x) - \int k(x,\cdot) d\tg(x)=\frac{1}{2}\nabla\MMD^2(\mu,\tgX)$. We have:	
\begin{align*}	
\ps{\nabla \MMD^2(\mu, \tg)& - \nabla \MMD^2(\nu, \tg), \mu-\nu} \le 	\|\nabla \MMD^2(\mu, \tg) -\nabla \MMD^2(\nu, \tg)\|_{\infty}\|\mu-\nu\|_{TV}\nonumber\\	
&\le 2\|f_{\mu,\tg}-f_{\nu,\tg}\|_{\infty}\|\mu-\nu\|_{TV} \\	
&\le 2\sup_{y \in \X} |\int k(x,y)d\mu(x)-\int k(x,y)d\nu(x)|\|\mu-\nu\|_{TV}\quad \\\label{eq9}\\
\intertext{since by the reproducing property and Cauchy-Schwarz inequality, $k(x,y)=\langle k(x,\cdot),k(y,\cdot) \rangle \le \|k(x,\cdot)\|_k\|k(y,\cdot)\|_k=\sqrt{k(x,x)k(y,y)}\le c_k$, and $y\mapsto k(x,y)$ is measurable,}
	&\le 2 c_k \sup_{\substack{f:\X \rightarrow[-1,1] \\f \text{mesurable}}} |\int f(x)d\mu(x)-\int f(x)d\nu(x)|\|\mu-\nu\|_{TV}\\
	&\le 2c_k \|\mu-\nu\|_{TV}^2 \le 4c_k(\KL(\mu|\nu)+\KL(\nu|\mu))=4c_k\ps{\nabla \phine(\mu)-\nabla \phine(\nu), \mu-\nu},	\end{align*}	where the last inequality results from Pinsker's inequality.	We conclude by using \Cref{lem:equivalences}.\end{proof}

\begin{remark}
(Case of neural network optimization). It is interesting to quantify the constant $c_k$ for some kernels of interest, for instance when optimizing an 
infinite-width one hidden layer neural network as in \citet{arbel2019maximum}. Consider a regression task where the labelled data $(z,y)\sim P$ where $P$ denotes some fixed data distribution. For any input $z$, the output of a single hidden layer neural network parametrized by $w \in \X$ can be written:
\begin{equation*}
    f_w(z)=\frac{1}{N}\sum_{j=1}^N a_j \sigma(\langle b_j, z\rangle) %
   = \int_{\X} \phi(z,w)d\mu(w), %
\end{equation*}
where $a_j$ and $b_j$ denote output and input weights of neuron $j=1,\dots,N$ respectively,  $w_j=(a_j,b_j)$ and $\mu = \nicefrac{1}{N}\sum_{j=1}^N \delta_{w_j}$. In the infinite-width setting, the limiting risk in this regression setting is written for any distribution $\mu\in \cP(\X)$ on the weights as $\E_{(z,y)\sim P}[\|y - \int \phi(z,w)d\mu(w)\|^2]$.   When the model is well-posed, i.e. there exists a distribution $\mu^*$ over weights such that $\mathbb{E}[y|z=\cdot]=\int\phi(\cdot,w)d\mu^*(w)$, then  the limiting risk writes as an MMD with $k(w,w')= \mathbb{E}_{z\sim P}[\phi(z,w)^T \phi(z,w)]$ (\citet[Proposition 20]{arbel2019maximum}). Hence, bounding $c_k = \sup_{w\in \X}k(w,w)$ depends on the choice of the activation function $\sigma$ and on bounding  the output weights. If $\sigma$ is bounded (e.g. $\sigma$ is the sigmoid activation) then bounding $c_k$ corresponds to bounding the output weights. If $\sigma$ is the RelU activation, then bounding $c_k$ depends on bounding both input and output weights as well, and on the data distribution $P$. 
\end{remark}

\section{Related work - Optimization over measures using the Wasserstein geometry} \label{sec:related_work}

In this section, we attempt to clarify the differences between the (Radon) vector space geometry considered in this paper and the Wasserstein geometry, developed in particular in \cite{otto2001geometry,villani2003topics,ambrosio2008gradient}.

Given an optimisation problem over $\mathcal{P}(\X)$ the set of probability distributions over $\X$, one can consider different geometries over $\mathcal{P}(\X)$. The one adopted in our paper casts $\mathcal{P}(\X)$ as a subset of a normed space of measures, such as $L^2(\rho)$ where $\rho$ is a reference measure, or Radon measures. In this space, the shortest distance paths between measures are given by their square-norm distance. Moreover in this setting, one can consider the duality of measures with continuous functions and the mirror descent algorithm, as we do in this work.

In contrast, another possibility is to restrict $\mathcal{P}(\X)$ to the probability distributions with bounded second moments, denoted $\mathcal{P}_2(\X)$, equipped with Wasserstein-2 ($W_2$) distance. The space $(\mathcal{P}_2(\X),W_2)$, called the Wasserstein space, is a metric space equipped with a rich Riemannian structure (often referred to as "Otto calculus") where the shortest distance paths are given by the $W_2$ distance and associated geodesics. In this setting, one can leverage the Riemannian structure to discretize ($W_2$) gradient flows and consider algorithms such as ($W_2$) gradient descent, in analogy with Riemannian gradient descent.

While both frameworks yield optimisation algorithms on measure spaces, the geometries and algorithms are very different. Both the notion of convexity  (along $L^2$ versus $W_2$ geodesics) and of gradients  (first variation vs gradient of first variation) differ; and by extension so do many definitions. Consequently, the conditions needed for the  convergence of mirror descent and $W_2$ gradient descent over an objective functional $\mathcal{F}$ greatly differ since they rely on the chosen geometry through the definitions of convexity, smoothness, or differentiability.

Wasserstein gradient descent should be thought of the analog of Riemannian gradient descent in infinite dimensions. Consequently, mirror descent yields updates on measures allowing for change of mass (see \Cref{rmk:dual_iteration}), while $W_2$ gradient descent preserves the mass, since the updates on measures write as pushforwards (i.e., displacement of particles supporting the measures). To summarize, the mirror descent scheme we consider is very different in nature to the gradient descent schemes based on the Wasserstein geometry (e.g. \cite{chizat2018global, mei2018mean,rotskoff2018neural,wibisono2018sampling, korba2020non,salim2020wasserstein, korba2021kernel}), due to the different geometry.

\section{Proofs}
\subsection{Proof of \Cref{th:rate}}\label{sec:proof_th_rate}

\begin{proof}
	Since $\cF$ is $L$-smooth relative to $\phi$ over $\cR$ and we assumed that $(\mu_n)_{n\in\N}\in \cR^\N$, we have
	\begin{equation}\label{eq:conv_smooth_new}
		\cF(\mu_{n+1})\le \cF(\mu_n) + d^+\!\cF(\mu_n)(\mu_{n+1}	-\mu_n)+L D_{\phi}(\mu_{n+1}|\mu_n).
	\end{equation}
	Applying Lemma \ref{lem:three-point} to the convex function $\cG_n(\nu)=\frac{1}{L}d^+\!\cF(\mu_n)(\nu-\mu_n)$, with $\mu=\mu_n$ and $\bar{\nu}=\mu_{n+1}$ yields
	\begin{equation*}
		d^+\!\cF(\mu_n)(\mu_{n+1}	-\mu_n) + L D_{\phi}(\mu_{n+1}|\mu_n) \le d^+\!\cF(\mu_n)(\nu-\mu_n) + LD_{\phi}(\nu| \mu_n) - L D_{\phi}(\nu| \mu_{n+1}).
	\end{equation*}
	Fix $\nu\in\cR$, then \eqref{eq:conv_smooth_new} becomes:
	\begin{equation}\label{eq:conv_3point_new}
		\cF(\mu_{n+1})\le \cF(\mu_n) + d^+\!\cF(\mu_n)(\nu-\mu_n)	+ LD_{\phi}(\nu| \mu_n) - L D_{\phi}(\nu| \mu_{n+1}).
	\end{equation}
	This shows in particular, by substituting $\nu=\mu_n$ and since $D_{\phi}(\nu| \mu_{n+1})\ge 0$, that $\cF(\mu_{n+1})\le \cF(\mu_n)$, i.e. $\cF$ is decreasing at each iteration. Since $\cF$ is $l$-strongly convex relative to $\phi$, we also have:
	\begin{equation}
		d^+\!\cF(\mu_n)(\nu-\mu_n)\le \cF(\nu)-\cF(\mu_n) - l D_{\phi}(\nu|\mu_n) 
	\end{equation}
	and \eqref{eq:conv_3point_new} becomes:
	\begin{equation}\label{eq:conv_strcvx_new}
		\cF(\mu_{n+1})\le \cF(\nu) + (L-l) D_{\phi}(\nu| \mu_{n}) - LD_{\phi}(\nu|\mu_{n+1}).
	\end{equation} 
	By induction, similarly to \citet{lu2018relatively}, we sum \eqref{eq:conv_strcvx_new} over $n$, obtaining
	\begin{align*}
	\sum_{i=1}^{n} \left(\frac{L}{L-l}\right)^i \cF(\mu_{i})&\le \sum_{i=1}^{n} \left(\frac{L}{L-l}\right)^i \cF(\nu) + L  D_{\phi}(\nu| \mu_{0}) - L \left(\frac{L}{L-l}\right)^n D_{\phi}(\nu|\mu_{n})
	\end{align*}
	Using the monotonicity of $(\cF(\mu_n))_{n\ge 0}$ and the positivity of $D_{\phi}(\nu|\mu_n)$, we have
	\begin{equation*}
	    \sum_{i=1}^{n} \left(\frac{L}{L-l}\right)^i \left(\cF(\mu_n) - \cF(\nu) \right)\le L  D_{\phi}(\nu| \mu_{0})- L \left(\frac{L}{L-l}\right)^n D_{\phi}(\nu|\mu_{n})\le L  D_{\phi}(\nu| \mu_{0}).\qedhere
	\end{equation*}
\end{proof}

\subsection{Proof of  \Cref{prop:sinkhorn_mirror_descent}}\label{sec:proof_prop_sinkhorn_mirror_descent}

Let $\pi,\bar{\pi} \in \cP(\X\times \Y)$, $h>0$ and $\xi=\bar{\pi}-\pi$ hence any integral with respect to $\xi$ of constant functions is null. We have:
\begin{align*}
 \Fsink(\pi + h \xi)& -\Fsink(\pi)  =\KL(\px(\pi+h\xi)|\tgX) -\KL(\px \pi|\tgX)  \\
& = \int \log \left(\px  \pi +h \px\xi\right)d\px(\pi + h\xi) - \int\log (\tgX) d\px(\pi+h\xi) - \int \log\left(\frac{\px\pi}{\tgX}\right)d\px\pi\\
 &= h \int \log\left(\frac{\px\pi}{\tgX}\right) d\px\xi +\underbrace{ \int \log \left(1+h \frac{\px\xi}{\px\pi}\right)d\px\pi}_{\approx h\int \frac{\px\xi}{\px\pi}d\px\pi + o(h) = 0+o(h)} +\underbrace{ h \int \log \left(1+h \frac{\px\xi}{\px\pi}\right)d\px\xi}_{\approx h^2 \int \frac{\px\xi}{\px\pi}d\px\xi+ o(h^2)}\\
& = h \int \log\left(\frac{\px\pi}{\tgX}\right) d\px\xi +o(h).
\end{align*}
Consequently,
\begin{equation}
   \lim_{h\to 0^+}\frac{ \Fsink(\pi + h \xi) -\Fsink(\pi)}{h}= \int_{\X} \log\left(\frac{\px\pi}{\tgX}\right)d\xi = \int_{\X \times \Y} \log\left(\frac{\px\pi}{\tgX}\right)d\xi.
\end{equation}
Hence, when it exists, $\nabla_{\!C}\Fsink(\pi)= \ln(d\px\pi/d\tgX)$. Moreover, the sets $\Pi(*,\tgY)$ and $\Pi(\tgX,*)$ are $L^\infty$-weak-* closed.\footnote{Indeed, take $(\pi_n)_{n\in N}\in \Pi(*,\tgY)$ converging weakly to some $\bar\pi$. As $\langle g,\pi_n\rangle_{\X\times\Y}=\langle g,\tgY\rangle_\Y$ for all $g\in L^\infty(\Y,\R)$, we obtain that $\langle g,\bar\pi\rangle_{\X\times\Y}=\langle g,\tgY\rangle_\Y$ which precisely means that $\py \bar\pi=\tgY$.} Besides, $\KL$ has weak-* compact sublevel sets and is weak-* l.s.c. Hence \citet[Theorem 3.2.2]{attouch2014variational} applies, and the iterates $(\pi_n)_{n\ge 0}$ exist, as originally shown by \citet{csiszar1975}. As $\pi_n=e^{(f+g-c)/\epsilon}\tgX\otimes\tgY$~\citep[Section 6]{Nutz2021IntroductionTE} with $f\in L^\infty(\X)$ and $g\in L^\infty(\Y)$, we have that $x\mapsto\ln(d\mu_n(x)/d\tg(x))\in L^\infty(\X,\R)$; indeed as $c\in L^{\infty}$, the first marginal $\mu_n$ of $\pi_n$ is an integral of functions bounded by strictly positive quantities.

Consider a coupling $\pi\in\cP(\X\times\Y)$ with $\pi\ll\pi_n$ and denote by $\mu$ its first marginal.  We have that $\Fsink(\pi_n)=\int_\X \ln(\mu_n/\tg)\,d\mu_n$ and $\ps{\nabla_{\!C} \Fsink(\pi_n),\pi-\pi_n}=\iint\ln(d\mu_n(x)/d\tg(x)) (\pi(dx,dy)-\pi_n(dx,dy))$. Simplifying and using \eqref{eq:disintegration} twice we obtain the identity:
\begin{align*}
        \Fsink&(\pi_n)+\ps{\nabla_{\!C} \Fsink(\pi_n),\pi-\pi_n}+\KL(\pi|\pi_{n})\\
        &=\int\ln(d\mu_n/d\tg(x)) \mu_n(dx)+\iint\ln(d\mu_n/d\tg(x)) \pi(dx,dy)\\
        &\hspace{1cm}-\iint\ln(d\mu_n/d\tg(x)) \pi_n(dx,dy)+\KL(\pi|\pi_{n})\\
        &=\int\ln(d\mu_n/d\tg(x)) \mu(dx)+\KL(\px \pi|\mu_n)+\KL(\pi|\px  \pi \otimes \nicefrac{\pi_n}{\mu_n})\\
        &=\int\ln(d\mu/d\tg(x)) \mu(dx)+\KL(\pi|\px  \pi \otimes \nicefrac{\pi_n}{\mu_n})=\KL(\pi|\tgX\otimes\nicefrac{\pi_n}{\mu_n})=\KL(\pi|\pi_{n+\frac12}).
    \end{align*}
   We conclude by taking the argmin over $\pi\in C$.

\subsection{Proof of \Cref{prop:quantitative-stability}}\label{sec:proof_quantitative-stability}

The proof of~\Cref{prop:quantitative-stability} essentially relies on bounding the entropic potentials by the marginals, as in~\citet{luise2019sinkhorn}[Theorem C.4]. For their purpose \citet{luise2019sinkhorn} assume that $\X=\Y$ and that $c(x,y)=c(y,x)$. These assumptions are not needed here so we revisit their proof and show that their bound holds for general bounded costs. We define 
\begin{equation*}
	D_c = \frac12 \sup [c(x,y)+c(x',y')-c(x,y')-c(x',y)],
\end{equation*}
where the supremum runs over $x,x'\in \X$ and $y,y'\in\Y$. 
When $\mu$ and $\nu$ are probability measures on $\X$ and $\Y$ respectively we define the soft-$c$ transform mappings $T_\mu\colon L^\infty(\X)\to L^\infty(\Y)$ and $T_\nu\colon L^\infty(\Y)\to L^\infty(\X)$ by 
\[
	T_\mu(f)(y)=-\eps\ln\Big(\int_\X e^{(f(x)-c(x,y))/\eps}\mu(dx)\Big)
\]
and 
\[
	T_\nu(g)(x)=-\eps\ln\Big(\int_\Y e^{(g(y)-c(x,y))/\eps}\nu(dy)\Big).
\]
These mappings arise naturally in the context of Sinkhorn's algorithm since if $\pi\in\PiexpInf$ has marginals $(\mu,\nu)$, we can write  $\pi(dx,dy)=e^{(f(x)+g(y)-c(x,y))/\eps}\mu(dx)\nu(dy)$ and taking marginals implies
\begin{equation} \label{eq:fixed-points-from-marginals}
	g=T_\mu(f)\quad\text{and}\quad f=T_\nu(g).
\end{equation}

\citet{luise2019sinkhorn} use the \emph{Hilbert metric} to prove their result, a classical tool to analyze matrix scaling problems~\citep{Franklin1989}, which for our purpose here reduces to the following semi-norm.

\begin{definition}
    When $f\in L^\infty(\X)$ we set $\normvar{f}=(\sup_\X f)-(\inf_\X f)$. We similarly define $\normvar{g}$ for $g\in L^\infty(\Y)$.
\end{definition}

We are now ready to state our version of~\citet{luise2019sinkhorn}[Theorem C.4].

\begin{proposition} \label{quantitative-stability-estimate}
	Let $\pi,\pit\in\PiexpInf$ with marginals $(\mu,\nu)$ and $(\mut,\nut)$ respectively. Write $\pi=e^{(f+g-c)/\eps}\mu\otimes\nu$ and $\pit=e^{(\ft+\gt-c)/\eps}\mut\otimes\nut$. Then 
	\[
		\normvar{f-\ft}+\normvar{g-\gt}\le 2\epsilon \,e^{3D_c/\epsilon}\big(\normtv{\mu-\mut} + \normtv{\nu-\nut}\big).
	\] 
\end{proposition}

The proof of this quantitative stability estimate mainly relies on the classical result that the soft $c$-transform mappings are contractions in the Hilbert metric; this result is at the heart of the proof of the classical linear convergence rate of Sinkhorn~\citep[see][for a proof]{Franklin1989,chen2016entropic}.

\begin{proposition} \label{lemma:contraction-Hilbert-metric}
    $\normvar{T_\mu(\ft)-T_\mu(f)}\le \lambda\normvar{\ft-f}$ with $\lambda=\frac{e^{D_c/\eps}-1}{e^{D_c/\eps}+1}<1$. 
\end{proposition}

We will also need the following lemma which is essentially contained in~\citet{luise2019sinkhorn}.

\begin{lemma} \label{lemma:logbound}
	Let $f=T_\nu(g)$ for some $g\in L^\infty(\Y)$. Then
	\[
		\normvar{T_\mut(f)-T_\mu(f)}\le 2\epsilon \,e^{2D_c/\epsilon}\normtv{\mu-\mut}.
	\]
	Likewise if $g=T_\mu(f)$ for some $f\in L^\infty(\X)$, 
	\[
		\normvar{T_\nut(g)-T_\nu(g)}\le 2\epsilon \,e^{2D_c/\epsilon}\normtv{\nu-\nut}.
	\]	
\end{lemma}
\begin{proof}[Proof of \Cref{lemma:logbound}]
	For any $f\in L^\infty(\X)$ we have by definition
	\begin{align*}
		T_\mut(f)(y)-T_\mu(f)(y) &= \epsilon\log\Big(\int_\X e^{(f(x)-c(x,y))/\epsilon}\mu(dx)\Big) - \epsilon\log\Big(\int_\X e^{(f(x)-c(x,y))/\epsilon}\mut(dx)\Big). 
	\end{align*}
	To control this difference of logs,~\citet{luise2019sinkhorn}[Lemma C.2] use the bound $\abs{\log(a)-\log(b)}\le \max\{a^{-1},b^{-1}\}\abs{a-b}$ (for any $a,b>0$). We have $\int_\X e^{(f(x)-c(x,y))/\epsilon}\mu(dx)\ge e^{\inf_x [f(x)-c(x,y)]/\epsilon}$ and the same lower bound holds for $\int_\X e^{(f(x)-c(x,y))/\epsilon}\mut(dx)$. Therefore 
	\begin{align*}
		\abs{T_\mut(f)(y)-T_\mu(f)(y)} &\le \epsilon\,e^{-\inf_x[f(x)-c(x,y)]/\epsilon} \int_\X e^{(f(x')-c(x',y))/\epsilon}\abs{\mu-\mut}(dx')\\
		&\le \epsilon\,e^{\sup_x[c(x,y)-f(x)]/\epsilon} e^{\sup_{x'}[f(x')-c(x',y)]/\epsilon} \normtv{\mu-\mut}.
	\end{align*}
	This implies when taking the supremum over $y\in\Y$
	\[
		\normvar{T_\mut(f)-T_\mu(f)}\le 2\norm{T_\mut(f)-T_\mu(f)}_\infty\le 2\epsilon\,e^{\sup_{x,x',y}[f(x')-f(x)+c(x,y)-c(x',y)]/\epsilon}\normtv{\mu-\mut}.
	\]
	This last inequality is valid for any $f\in L^\infty(\X)$. If in addition we take $f$ to be an image $f=T_\nu(g)$, then we have the standard estimate for any given $x,x' \in\X$ and $y\in\Y$
	\begin{align*}
		-f(x)&=\eps\ln\Big(\int_\Y e^{(g(y')-c(x,y'))/\eps}\nu(dy')\Big)\\
		&=\eps\ln\Big(\int_\Y e^{(c(x,y)+c(x',y')-c(x',y)-c(x,y'))/\eps}e^{(g(y')-c(x',y'))/\eps}\nu(dy')\Big) +c(x',y)-c(x,y)\\
		&\le 2D_c + \eps\ln\Big(\int_\Y  e^{(g(y')-c(x',y'))/\eps}\nu(dy')\Big)+c(x',y)-c(x,y)\\
		&= 2D_c - f(x')+c(x',y)-c(x,y).
	\end{align*}
	This shows that  $\sup_{x,x',y}[f(x')-f(x)+c(x,y)-c(x',y)] \le 2D_c$. As a consequence,
	\[
		\normvar{T_\mut(f)-T_\mu(f)}\le 2\epsilon\,e^{2D_c/\epsilon}\normtv{\mu-\mut}.
	\]

	By symmetry the corresponding bound can be derived for quantities on $\Y$. 
\end{proof}

\begin{proof}[Proof of~\Cref{quantitative-stability-estimate}]
	Having in mind the fixed point equations~\eqref{eq:fixed-points-from-marginals} for $(f,g)$ and the corresponding ones for $(\ft,\gt)$ we write 
    \begin{align*}
        \normvar{\ft-f} &= \normvar{T_\nut(\gt)-T_\nu(g)}\\
        &\le \normvar{T_\nut(\gt)-T_\nut(g)} + \normvar{T_\nut(g)-T_\nu(g)},
    \end{align*}
	and similarly, $\normvar{\gt-g}\le \normvar{T_\mut(\ft)-T_\mut(f)} + \normvar{T_\mut(f)-T_\mu(f)}$.
    By~\Cref{lemma:contraction-Hilbert-metric}, $\normvar{T_\mut(\ft)-T_\mut(f)}\le\lambda\normvar{\ft-f}$ and $\normvar{T_\nut(\gt)-T_\nut(g)}\le\lambda\normvar{\gt-g}$. Combining, we obtain
	\begin{equation*}
		(1-\lambda) \big(\normvar{\ft-f} + \normvar{\gt-g}\big) \le \normvar{T_\mut(f)-T_\mu(f)} + \normvar{T_\nut(g)-T_\nu(g)}.
	\end{equation*}
	\Cref{lemma:logbound} takes care of the right-hand side, and this results in
	\begin{equation*}
		(1-\lambda) \big(\normvar{\ft-f} + \normvar{\gt-g}\big) \le 2\epsilon \,e^{2D_c/\epsilon}\big(\normtv{\mu-\mut} + \normtv{\nu-\nut}\big).
	\end{equation*}
	Finally we divide by $1-\lambda$ and bound $(1-\lambda)^{-1}=(e^{D_c/\eps}+1)/2\le e^{D_c/\eps}$.
\end{proof}

\begin{proof}[Proof of~\Cref{prop:quantitative-stability}]
    Let $\pi,\pit\in\PiexpInf$ with marginals $(\mu,\tgY)$ and $(\mut,\tgY)$ respectively. Write $\pi=e^{(f+g-c)/\eps}\mu\otimes\tgY$ and $\pit=e^{(\ft+\gt-c)/\eps}\mut\otimes\tgY$. We emphasize that $\pi$ and $\pit$ have the same second marginal $\tgY$. Then
    \begin{align*}
        &\eps\KL(\pit|\pi) = \eps\KL(\pit|\pi) +\eps\KL(\pi|\pit) - \eps\KL(\pi|\pit) \\
        &= \iint (\tilde f-f+\tilde g-g+\eps\ln\Big(\frac{d\mut}{d\mu}\Big))\,d\pit+\iint (f-\ft+g-\gt+\eps\ln\Big(\frac{d\mu}{d\mut}\Big))d\pi -\eps\KL(\pi|\pit)\\
        &= \iint (\tilde f-f+\tilde g-g)\,(d\pit-d\pi) +\eps\KL(\mut|\mu)+\eps\KL(\mu|\mut)-\eps\KL(\pi|\pit).
    \end{align*}
   Part of the first term vanishes since $\iint(\gt-g)\,(d\pit-d\pi)=\int_\Y(\gt-g)\,(d\tgY-d\tgY)=0$, and we can get rid of the last two terms by using the data processing inequality $\KL(\mu|\mut)\le\KL(\pi|\pit)$. Thus 
    \[
        \eps\KL(\pit|\pi)\le \normvar{\tilde f-f}\normtv{\mut-\mu} +\eps\KL(\mut|\mu).
    \]
    Applying~\Cref{quantitative-stability-estimate} we obtain
    \[
        \eps\KL(\pit|\pi)\le 2\eps e^{3D_c/\epsilon}\normtv{\mut-\mu}^2 +\eps\KL(\mut|\mu),
    \]
    and after dividing by $\eps$, Pinsker's inequality yields 
    \begin{equation*}
        \KL(\pit|\pi)\le (1+4 e^{3D_c/\epsilon})\KL(\mut|\mu).\qedhere
    \end{equation*}
    
\end{proof}

\subsection{Proof of  \Cref{prop:em-mirror-descent}}\label{sec:proof_enveloppe_FEM}

\begin{proposition}[EM as mirror descent]\label{prop:em-mirror-descent-precise}Let $C=\Pi(*,\tgY)$. Assume that for all $\pi \in C$ there exists a $q_*(\pi)\in \Q$ solving \eqref{eq:def-FEM}, that, for $p_h=p_{q_*((1-h)\pi_n+h\pi)}$, $\nicefrac{d\pi_n}{dp_h}$ converges pointwise to $\nicefrac{d\pi_n}{dp_{q_n}}$ for $h\rightarrow 0^+$ with $|\ln(\nicefrac{d\pi_n}{dp_{h}})|\le \cG_n$ for some $ \cG_n\in L^1(\pi+\pi_n)$, that $\nicefrac{d\pi_n}{dp_{q_n}}(\cdot,\cdot)\in[a_n,b_n]$ for some $a_n>0$ and $b_n>0$, and that $\sup_{q\in \Q}|\ln(\nicefrac{d\pi_n}{d p_q})|<\infty$. Then the EM iterations~\eqref{eq:em-iterations-M-step}--\eqref{eq:em-iterations-E-step} can be written as a mirror descent iteration with objective function $\FEM$, Bregman potential $\phine$ and constraints $C$,

    \begin{equation}\label{eq:em_md_iterate-precise}
        \pi_{n+1}=\argmin_{\pi\in C} \ps{\nabla_{\!C} \FEM(\pi_n),\pi-\pi_n} + \KL(\pi|\pi_n),
        \end{equation}
    
         with $\nabla_{\!C}\FEM(\pi_n)= \ln(d\pi_n/dp_{\qEM_n})\in L^\infty(\X\times \Y)$.
\end{proposition}

\begin{remark}
Note that our assumptions on the sequence $(p_h)_{h\in[0,1]}$ in \Cref{prop:em-mirror-descent-precise} are very similar to what the fundamental theorem of $\Gamma$-convergence would provide (see \cite{dal1987gamma}, \citet[Theorem 2.10]{braides2002gamma}). It is indeed straightforward to prove $\Gamma$-convergence \citep[see][Theorem 2.1]{braides2002gamma} of the sequence $(f_{n,\pi}(\cdot,h))_{h\in[0,1]}$ in $h=0^+$ with $f_{n,\pi}(p,h):=\KL(\pi_n+h(\pi-\pi_n)|p)$, owing to the convexity and joint weak-* lower semicontinuity of $\KL$. However, to prove the convergence of the sequence of minimizers $(p_h)_{h\in[0,1]}$, one would need the equicoercivity of $(\KL(\pi_h | p))_{h\in[0,1]}$ over $p\in \cP_Q$ \citet[Definition 2.9]{braides2002gamma}), which heavily depends on the properties of $\cP_Q$, e.g.\ considering a weak-* compact $\cP_Q$ would entail equicoercivity. 
\end{remark}

\begin{proof}
We will use here the envelope theorem to differentiate $\FEM$ and compute its first variation. 
We are going to apply \citet[Theorem 3]{Milgrom02envelopetheorems} leveraging properties of $\KL$. \citet[Theorem 3]{Milgrom02envelopetheorems} is written for the set $[0,1]\times X$, where $X$ is some set optimized over. Here $X=\cP_Q:=\{p_q \, | \, q\in\Q\}$ and the interval $[0,1]$ will be merely the scalar of the directional derivative we consider.

Let $n\ge 0$ and $\pi\in \cP(\X\times\Y)$. For $h\in[0,1]$, set $f_{n,\pi}(p,h):=\KL(\pi_n+h(\pi-\pi_n)|p)$ and $V_{n,\pi}(h)=\inf_{p\in \cP_Q} \KL(\pi_n+h(\pi-\pi_n)|p)$ to match the notations of \citet[Theorem 3]{Milgrom02envelopetheorems}. We have to show some equidifferentiability over $q\in \cP_Q$. Notice that the following expression does not depend on $p$,
\begin{align*}
    &\frac{1}{h}\left[\iint \ln\left(\frac{\pi_n+h(\pi-\pi_n)}{p}\right)d(\pi_n+h(\pi-\pi_n))-\iint \ln\left(\frac{\pi_n}{p}\right)d \pi_n\right]-\iint \ln\left(\frac{\pi_n}{p}\right)d(\pi-\pi_n)\\
    &=\frac{1}{h}\iint \ln\left(1+h\frac{(\pi-\pi_n)}{\pi_n}\right)d \pi_n=\frac{1}{h}[h\iint d(\pi-\pi_n)+O(h^2)]=0+O(h),
\end{align*}
so that we do have equidifferentiability when $h\to 0^+$. Our assumptions then allow to apply \citet[Theorem 3]{Milgrom02envelopetheorems}. We thus obtain that
\begin{equation*}
    d^+\FEM(\pi_n)(\pi-\pi_n)=d^+V_{n,\pi}(0)=\lim_{h\rightarrow 0^+}\iint \ln\left(\frac{d\pi_n}{dp_h}\right)d \pi_n.
\end{equation*}
Since $|\ln(\nicefrac{d\pi_n}{dp_{h}})|\le \G_n\in L^1(\pi+\pi_n)$ and $\nicefrac{d\pi_n}{dp_h}$ converges pointwise to $\nicefrac{d\pi_n}{dp_{q_n}}$ for $h\rightarrow 0^+$ (recall that $q_n= q_*(\pi_n)$ by definition), we can apply the dominated convergence theorem to interchange the limit and the integral. Consequently $d^+\FEM(\pi_n)(\pi-\pi_n)=\iint \ln(\nicefrac{d\pi_n}{dp_{q_n}})d(\pi-\pi_n)$ proving that $\nabla \FEM(\pi_n)=\ln(\nicefrac{d\pi_n}{dp_{q_n}})\in L^\infty$ since $\nicefrac{\pi_n}{p_{q_n}}(\cdot,\cdot)\in[a_n,b_n]$ for some $a_n>0$ and $b_n>0$.

Then, for $\pi_n$ the EM iterate at time $n$, and for any coupling $\pi$, we have the identity:
\begin{align*}
  &\FEM(\pi_n)+\ps{\nabla_{\!C}\FEM(\pi_n),\pi-\pi_n}+\KL(\pi|\pi_n) \\
  &= \int \ln(d\pi_n/dp_{\qEM_n}(x))\pi_n(dx) + \int \ln(d\pi_n/dp_{\qEM_n}(x))\, (\pi-\pi_n)(dx) + \int \ln(d\pi/d\pi_n(x))\pi(dx) \\
  &= \int \ln(d\pi/dp_{\qEM_n}(x))\pi(dx) = \KL(\pi|p_{\qEM_n}).
\end{align*}
Note that $\qEM_n$ is optimal in~\eqref{eq:def-FEM}, whence \eqref{eq:em_md_iterate} matches~\eqref{eq:em-iterations-E-step}.
\end{proof}

\subsection{Proof of  \Cref{prop:em_mirror_descent}}\label{sec:proof_first_variations_FEMK}

\begin{proof}
Let $\pi,\bar{\pi} \in \cP(\X\times \Y)$, $h>0$ and $\xi=\bar{\pi}-\pi$, so $\iint_{\X \times \Y} \xi(dx,dy)=0$. We have:
\begin{align*}
    &\FEMK(\pi+h\xi) - \FEMK(\pi) =   \KL(\pi+h\xi|\px(\pi+h\xi) \otimes K ) - \KL(\pi|\px\pi \otimes K )  \\
    &= \int \log\left( \frac{\pi +h \xi}{\px (\pi+h\xi) \otimes K} \right)d(\pi +h\xi) -   \int \log\left( \frac{\pi }{\px \pi \otimes K} \right)d\pi\\
   & = \int \log\left( \pi +h \xi  \right)d(\pi +h\xi) - \int \log\left( \px(\pi +h \xi) \otimes K \right)d(\pi +h\xi)
    -   \int \log\left( \pi \right)d\pi + \int \log \left(\px \pi \otimes K \right)d\pi \\
   & = h \int \log \pi d\xi + \int \underbrace{\log \left( 1+h \frac{\xi}{\pi} \right) d\pi}_{\approx h \int \frac{\xi}{\pi}d\pi + o(h) =0+ o(h)} + h \underbrace{\int \log \left( 1+h \frac{\xi}{\pi} \right) d\xi}_{\approx h^2 \int \frac{\xi}{\pi}d\xi +o(h^2)}\\
    & - h \int \log(\px \pi \otimes K)d\xi- \underbrace{\int \log \left( 1+h \frac{\px \xi \otimes K}{\px \pi \otimes K} \right)d\pi}_{\approx h \int \frac{\px \xi \otimes K}{\px \pi \otimes K} d\pi+o(h) = h\int \frac{\px \xi}{\px \pi} d\pi +o(h)= 0+o(h)} -\underbrace{ h \int \log\left( 1+h \frac{\px \xi \otimes K}{\px \pi \otimes K} \right) d\xi}_{h^2\int \frac{\px \xi \otimes K}{\px \pi \otimes K} d\xi+o(h^2) }\\ %
   & = h \int \log\left(\frac{\pi}{\px\pi \otimes K}\right)d\xi+o(h).
\end{align*}
Hence
\begin{align*}
   \lim_{h\to 0^+} \frac{\FEMK(\pi+h\xi) - \FEMK(\pi) }{h} = \int \log\left(\frac{\pi}{\px\pi \otimes K}\right)d\xi. %
\end{align*}

To show that $\nabla \FEMK(\pi_n)$ belongs to $L^\infty$, we proceed by induction. Let $n\ge 0$ and assume that $T_K\mu_n\gg\tgY$ and $\mu_n=e^{f_n(x)}\tgX$ with $f_n$ bounded (which we explicitly assumed for $\mu_0$) then the multiplicative update \eqref{eq:em-latent} shows that $f_{n+1}$ has the same property. Furthermore \eqref{eq:em-latent} gives
\begin{align*}
    \frac{\pi_{n+1}}{\mu_{n+1}\otimes K}(\cdot)=\frac{\mu_{n}(\cdot) \frac{K(\cdot,dy)\tgY(dy)}{\int_{\mathcal{X}} K(x,dy)\mu_{n}(dx)}}{\mu_{n}(\cdot)\otimes K(\cdot,dy)
    \int_{\mathcal{Y}} \frac{K(\cdot,dy')\tgY(dy')}{\int_{\mathcal{X}} K(x,dy')\mu_{n}(dx)}}=\frac{\tgY(dy)}{\int_{\mathcal{X}} K(x,dy)\mu_{n}(dx)\times \int_{\mathcal{Y}} \frac{K(\cdot,dy')\tgY(dy')}{\int_{\mathcal{X}} K(x,dy')\mu_{n}(dx)}}
\end{align*}
Since $K(x,dy)=e^{-c(x,y)}\tgY$ with $c$ uniformly bounded, $\frac{\pi_{n+1}}{\mu_{n+1}\otimes K}(\cdot)$ is also bounded above and below by positive constants (depending on $n$).
\end{proof}

\subsection{Proof of \Cref{prop:rates_em}}\label{sec:proof_rates_em}

\begin{proof}
By the disintegration formula \eqref{eq:disintegration}, 
\begin{equation}\label{eq:decompose_FEMK-app}
    \FEMK(\pi) = \KL(\tgY|\py (\px \pi \otimes K)) +\int \KL(\pi/\tgY| (\px \pi \otimes K )/ \py (\px \pi \otimes K) ) d\tgY
\end{equation}
 Let    $\pi_*=\mu_*(dx)k(x,dy)\tgY(dy)/(T_K\mu_*)(dy)$. First, for any $\pi \in \cP(\X\times \Y)$, we have $\py(\px \pi \otimes K)= \int \px\pi(dx)k(x,\cdot)=T_K (\px\pi)$, hence by definition of $\mu^*\in\argmin_{\mu}\KL(\tgY|T_k\mu)$, $\pi_*$ minimizes the first term in \eqref{eq:decompose_FEMK-app}. Second, this choice leads to $\pi_*/\tgY = \mu_*\otimes K/T_K \mu^*$, cancelling the nonnegative second term in \eqref{eq:decompose_FEMK-app}. Hence $\pi_*$ is a minimizer of $\FEMK$ and $\FEMK(\pi_*) = \KL(\tgY|T_K \mu_*)$. Moreover, $\KL(\cdot|\cdot)$ is convex in both arguments, and $\pi\mapsto\px \pi \otimes K$ is linear. Consequently the composition $\FEMK$ is convex in $\pi$ and so is $\Fsink$ by the same arguments (see also \Cref{prop:rel-smooth-Sink}). By \eqref{eq:disintegration} and linearity of the Bregman divergence, $\KL(\pi|\tilde{\pi})= D_{\Fsink}(\pi|\tilde{\pi}) + D_{\FEMK}(\pi|\tilde{\pi})$, hence $\FEMK$ is 1-relatively smooth w.r.t. $\phine$. Hence, \Cref{th:rate} yields:
  \[
        \FEMK(\pi_n)\le \FEMK(\pi_*)+\frac{\KL(\pi_*|\pi_0)}{n}.
    \]
Since $\pi_0\in \Pi(*,\tgY)$, $\pi_0=\mu_0(dx)k(x,dy)\tgY(dy)/(T_K\mu_*)(dy)$, 
 \begin{align*}
     &\KL(\pi_*|\pi_0)= \KL(\mu_*|\mu_0)+\iint \ln\left(\frac{k(x,dy)\tgY(dy)/(T_K\mu_*)(dy)}{k(x,dy)\tgY(dy)/(T_K\mu_0)(dy)}\right)\pi_*(dx,dy) \\
     &= \KL(\mu_*|\mu_0)+\int_\Y \ln\left(\frac{\tgY(dy)/(T_K\mu_*)(dy)}{\tgY(dy)/(T_K\mu_0)(dy)}\right)\tgY(dy) =  \KL(\mu_*|\mu_0)+\KL(\tgY|T_K\mu_*)-\KL(\tgY|T_K\mu_0).
 \end{align*}
 Finally, we use the inequality $$\KL(\tgY|T_K \mu_n)=\KL(\pi_n|\py(\px\pi_n \otimes K))\le \KL(\py\pi_n|\px\pi_n \otimes K)= \FEMK(\pi_n).$$
\end{proof}

\end{document}